\documentclass[a4paper,12pt]{article}
\usepackage{amsmath,amssymb} \usepackage{euler,eucal}
\usepackage[utf8]{inputenc}
\usepackage[T1]{fontenc}
\usepackage[english]{babel}
\usepackage{tgpagella}
\usepackage[unicode=true,plainpages=false,breaklinks=true]{hyperref}
\usepackage{graphicx}
\hypersetup{colorlinks=false,pdfborder={0 0 0.1},linkcolor=magenta,anchorcolor=magenta,urlcolor=blue,citecolor=blue,
          pdftitle={Parallel cross interpolation for high-precision calculation of high-dimensional integrals},
          pdfauthor={S. V. Dolgov, D. V. Savostyanov},
          pdfkeywords={high-dimensional problems, tensor train format, cross interpolation, Ising integrals},
}

\def\email#1{\texttt{#1}}
\usepackage[super]{nth}
\usepackage[left=20mm,right=20mm,top=25mm,bottom=25mm,a4paper]{geometry}
\usepackage{amsthm} \usepackage{mathtools}
\theoremstyle{definition}
\newtheorem{remark}{Remark}
\newtheorem{theorem}{Theorem}

\usepackage[ruled]{algorithm}
\usepackage{algorithmicx}
\usepackage{algpseudocode}
\algrenewcommand{\algorithmiccomment}[1]{\hfill\textcolor{black!60!blue}{\textsl{\% #1}}}

\usepackage{tikz} \usetikzlibrary{positioning,shapes,calc,3d}
\usepackage{pgfplots}
\input{./ising.sty}
\let\eps=\varepsilon 

\let\le=\leqslant 
\let\leq=\leqslant \let\geq=\geqslant

\def\d{\mathrm{d}}  \def\dx{\mathrm{d}x} 
\def\R{\mathbb{R}}
\def\I{\mathcal{I}} \def\J{\mathcal{J}} 
\def\II{\mathbb{I}} \def\JJ{\mathbb{J}}
\def\Int{I} \def\C{C} \def\D{D} \def\E{E} 
\def\O{\mathcal{O}}
\def\Rnd{\mathcal{L}}

\def\rank{\mathop{\mathrm{rank}}\nolimits}
\def\det{\mathop{\mathrm{det}}\nolimits}
\def\vol{\mathop{\mathrm{vol}}\nolimits}
\def\maxvol{\mathop{\mathrm{maxvol}}\nolimits}

\def\new{\star}
\def\spin{\sigma} \def\spinup{\uparrow} \def\spindown{\downarrow}
\def\x{\mathbf{x}}
\def\qmc{\mathbf{q}}
\def\shift{\mathbf{s}}
\def\neval{N_{\mathrm{eval}}}
\def\e{\left._{\times 10}\right.}

\allowdisplaybreaks
\begin{document}

\title{Parallel cross interpolation for high--precision calculation of high--dimensional integrals\thanks{
 This work is supported by EPSRC grants EP/M019004/1 (S.D.) and EP/P033954/1 (D.S.).
 D.S. was additionally supported by the Rising Stars grant at the University of Brighton.
}}
\author{Sergey Dolgov\footnotemark[2] and Dmitry Savostyanov\footnotemark[3]}
\date{27 March 2019}
\maketitle

\renewcommand{\thefootnote}{\fnsymbol{footnote}}
\footnotetext[2]{University of Bath, Claverton Down, BA2 7AY Bath, UK,        \email{s.dolgov@bath.ac.uk}}
\footnotetext[3]{University of Brighton, Lewes Road, BN2 4GJ, Brighton UK    \email{d.savostyanov@brighton.ac.uk}}

\begin{abstract}
 We propose a parallel version of the cross interpolation algorithm and apply it to calculate high--dimensional integrals motivated by Ising model in quantum physics.
 In contrast to mainstream approaches, such as Monte Carlo and quasi Monte Carlo, the samples calculated by our algorithm are neither random nor form a regular lattice.
 Instead we calculate the given function along individual dimensions (modes) and use this data to reconstruct its behaviour in the whole domain.
 The positions of the calculated univariate fibers are chosen adaptively for the given function.
 The required evaluations can be executed in parallel both along each mode (variable) and over all modes.
 
 To demonstrate the efficiency of the proposed method, we apply it to compute high--dimensional Ising susceptibility integrals, arising from asymptotic expansions for the spontaneous magnetisation in two--dimensional Ising mo\-del of ferromagnetism.
 We observe strong superlinear convergence of the proposed method, while the MC and qMC algorithms converge sublinearly.
 Using multiple precision arithmetic, we also observed exponential convergence of the proposed algorithm.
 Combining high--order convergence, almost perfect scalability up to hundreds of processes, and the same flexibility as MC and qMC, the proposed algorithm can be a new method of choice for problems involving high--dimensional integration, e.g. in statistics, probability, and quantum physics.

\end{abstract}

\section{Introduction} 
 High--dimensional integrals occur often in statistics and probability (in e.g. expectations with multivariate probability distributions \cite{mcmc-handbook-2011}, inverse problems with uncertainty \cite{stuart-bayes-2010} and many more)
  or quantum mechanics \cite{meyer-mctdh-1990}.
 Analytical formulae for them are rarely available, hence numerical approaches become the mainstream approach.
 Unfortunately, high--dimensional integrals are notoriously difficult for numerical methods as well.
 A na\"ive approach, based on tensor product of one--dimensional quadrature rules, requires the total number of function evaluations $N$ that grows exponentially with problem dimension $d,$ exceeding the possibilities of modern computers for $d\gtrsim 10.$
 This behaviour, known as the \emph{curse of dimensionality}, motivates development of special methods for the integration in higher dimensions.
 Currently the most popular methods are
  the Monte Carlo quadrature~\cite{metropolis-mc-1949},
  quasi Monte Carlo~\cite{nieder-qmc-1978,morokoff-qmc-1995,graham-QMC-2011,Schwab-HOQMC-2014},
  Markov chain Monte Carlo~\cite{mcmc-handbook-2011},  
  and their derivatives such as multilevel Monte Carlo methods~\cite{BarthSchwab-ml_mc_spde-2011,Schwab-MLQMC-2015,Scheichl-mlqmc-lognorm-2017,nobile-sg-mc-2016}.
 These algorithms are rigorously studied and many theoretical results are available, including error bounds which typically do not depend on problem dimension $d$ for problems of interest.
 Unfortunately, MC and qMC methods converge slowly --- the relative accuracy $\eps$ depends on the number of function evaluations $\neval$ as $\eps\sim \neval^{-\gamma},$ where the convergence rate $\gamma=0.5$ for MC and $0.5\leq\gamma\leq 1$ for qMC.
 The numerical costs therefore grow quickly when higher precision is required, making calculations expensive, prohibitively long, or impossible.
 Methods based on Smolyak's sparse grids~\cite{smolyak-1963,griebel-sparsegrids-2004,Bieri2009} are often used to mitigate, but can not fully remove, the curse of dimensionality.

 In this paper we consider a problem of numerical integration of a multivariate function in a simple tensor--product domain such as free space $\R^d$ or hypercube $[0,1]^d.$
 We follow the na\"ive approach and use a tensor product of univariate quadrature rules, hence reducing the problem to calculation and summation over the entries of a multi--dimensional array (which we call \emph{tensor}).
 To overcome the curse of dimensionality, we approximate the whole array based on a few entries from it, but avoid calculating the whole array.
 To achieve this, we develop and use the parallel version of the tensor cross interpolation algorithm proposed by one of the authors in~\cite{sav-qott-2014}.
 This algorithm interpolates the given array in the \emph{tensor train} (TT) decomposition~\cite{osel-tt-2011,ot-ttcross-2010}, essentially performing \emph{separation of variables}.
 The array entries are evaluated along one--dimensional lines or \emph{fibers}, each of which is formed by freezing all indices of the multivariate function and only varying one.
 The lines intersect forming \emph{crosses}, and on the positions of each cross the constructed approximation \emph{interpolates} the data exactly, which explains the name of the algorithm.
 The positions of the crosses, and hence the nodes of the quadrature rule, are chosen adaptively for the given function, following the maximum--volume method~\cite{gtz-maxvol-1997,gostz-maxvol-2010}.
 When the approximation is available, various observables, including the integral, can be computed in linear in $d$ time.

 Essentially, the proposed algorithm reconstructs all $n^d$ values of the function $f(x_1,\ldots,x_d)$ on a tensor product $n\times \cdots \times n$ quadrature grid from a linear in $d$ number of samples, which are adapted specifically to $f.$
 This adaptivity allows the proposed algorithm to locate important samples (e.g. areas of concentration of the density) and reach faster convergence, compared to mainstream numerical methods, such as MC and qMC, where the positions of the samples are either not optimised, or are optimal for a wide class of functions.
 For the family of Ising integrals, considered in the numerical experiments section of this paper, the proposed algorithm demonstrates high--order convergence of the order of $\eps\sim \neval^{-7},$ clearly outperforming MC and qMC. 
 Using multiple precision arithmetic, we were able to compute an integral in more than thousands dimensions to more than hundred decimal digits, observing  exponential convergence of the proposed method.
 As a flexible and non--intrusive algorithm, it can become a new method of choice for problems involving numerical integration in higher dimensions.

 Data-sparse algorithms based on tensor product decompositions 
   (canonical polyadic~\cite{sav-rr-2009}, 
   Tucker~\cite{ost-tucker-2008}, 
   tensor train (TT)~\cite{osel-tt-2011} 
   or Hierarchical Tucker (HT)~\cite{hk-ht-2009}) 
  have a long history of development~\cite{kolda-review-2009,hackbusch-2012,bokh-surv-2015,lars-review-2014},
  with applications in 
    quantum physics and chemistry~\cite{fkst-chem-2008,ost-chem-2010,sdwk-nmr-2014,dkh-cme-2014,ds-dmrgamen-2015},
    signal processing~\cite{dks-ttfft-2012,sav-rank1-2012},
    plasma modelling~\cite{dst-fb-2014},
    stochastics and uncertainty quantification~\cite{zyokd-ttcircuit-2015,dklm-tt-pce-2015},
    and fractional calculus~\cite{rst-volterra-2014,dpss-frac-2015}.
 However, scalable \emph{high performance} implementation of tensor algorithms is a relatively new area of research.
 A straightforward idea is to parallelise dense tensor algebra in computations of factors of a decomposition~\cite{solomonik-cyclops-2014}.
 However, this typically requires all-to-all communications which quickly limit scalability of MPI code.
 Another strategy is to parallelise a tensor decomposition over different factors, or \emph{dimensions}.
 One of the first examples of the latter was the parallel density matrix renormalization group (DMRG) algorithm~\cite{white-parallel-dmrg-2013} for ground state computations in quantum physics.
 In mathematical community this research direction started with dimension--parallel linear solver~\cite{etter-par-als-2016} and cross algorithms in HT format~\cite{grasedyck-par-cross-2015}.
 The main difficulty of parallelisation over dimension is the need of \emph{algorithmic modifications}, since state of the arts tensor algorithms were designed in intrinsically sequential way.
 Ideally, such modifications should not compromise numerical stability or convergence for the sake of parallel efficiency.

 In this paper we develop a parallel version of the TT cross interpolation algorithm~\cite{sav-qott-2014}.
 The parallel algorithm is adaptive and converges with the same rate as the sequential version, but involves only local communications with constant loading of processes, and demonstrates almost perfect scaling up to the ultimate partitioning where each process is responsible for a single direction (mode, variable).
 Moreover, further speedup can be achieved using OpenMP parallelisation of tensor algebra in each process.

 The rest of the paper is organised as follows.
 In Sec.~2 we recall the cross interpolation method for matrices and provide necessary definitions.
 In Sec.~3 we discuss how the matrix interpolation can be applied for high--dimensional arrays (tensors). 
  We compare currently existing methods and explain why the cross interpolation algorithm proposed by one of the authors in~\cite{sav-qott-2014} seems to be the most suitable for parallelisation over the dimensions.
  We then present the parallel version of this algorithm.
 In Sec.~4 we explain how the cross interpolation algorithm can be applied for numerical integration.
  We also introduce more formally the MC and qMC methods for the same purpose.
 In Sec.~5 we introduce Ising susceptibility integrals which will be our main example in this paper.
  We demonstrate that the proposed method achieves high--order (sometimes exponential) convergence, while the convergence of MC and qMC remains sublinear. 
  In the conclusion, we briefly summarise the results of this paper and discuss some challenges and potential directions for the future work. 
 
\section{Cross interpolation: notation, definitions and algorithms}
 \subsection{Cross interpolation of matrices}
  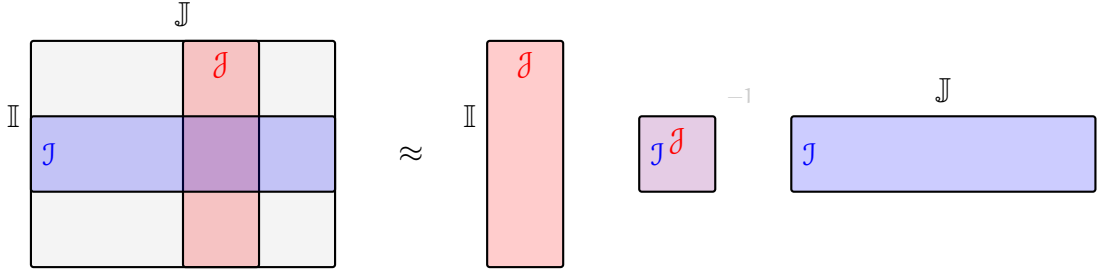
\begin{figure}[t]
   \begin{center}
    \begin{tikzpicture}[x=10mm,y=10mm]
     \tikzstyle{box}=[thick,rounded corners=1,draw=black,text opacity=1,fill=blue,fill opacity=.2];
     \draw[box,fill=black!20] (0,0) rectangle (4,3);
     \draw[box,fill=red]      (2,0) rectangle (3,3);
     \draw[box,fill=blue]     (0,1) rectangle (4,2);
     \node at (5,1.5){$\approx$};
     \draw[box,fill=red]      (6,0) rectangle (7,3);
     \draw[box,fill=red!50!blue]  (8,1) rectangle (9,2) node[above right]{$\scriptstyle -1$};
     \draw[box,fill=blue]     (10,1) rectangle (14,2);
     \node[left]           at  (0,2){$\II$};
     \node[left]           at  (6,2){$\II$};
     \node[above]          at  (2,3){$\JJ$};
     \node[above]          at (12,2){$\JJ$};
     \node[right,blue]     at  (0,1.5){$\I$};
     \node[right,blue]     at  (8,1.5){$\I$};
     \node[right,blue]     at (10,1.5){$\I$};
     \node[below,red]      at (2.5,3){$\J$};
     \node[below,red]      at (6.5,3){$\J$};
     \node[below,red]      at (8.5,2){$\J$};
    \end{tikzpicture}
   \end{center}
   \caption{Cross interpolation for matrices. 
   The full matrix $A$ is approximated by a low--rank decomposition $\tilde A$ based on a small number of columns and rows computed in $A.$
   Note that the approximation~\eqref{eq:mat} is exact in the positions of computed rows and columns.
   }
   \label{fig:mat}
  \end{figure} 
 Cross interpolation is based on a simple observation:
   for a given $m \times n$ matrix $A=[A(i,j)]_{i,j=1}^{m,n}$ its rank--$r$ interpolation can be recovered from 
   its $r$ columns $\J=\{\J^{(t)}\}_{t=1}^r$ and $r$ rows $\I=\{\I^{(s)}\}_{s=1}^r$ as follows:
 \begin{equation}\label{eq:mat}
  \begin{split}
  A(i,j)  \approx \tilde A(i,j) 
         & = \sum_{s=1}^r \sum_{t=1}^r A(i,\J^{(t)}) [A(\I,\J)]_{t,s}^{-1} A(\I^{(s)},j)
      \\ & = A(i,\J) [A(\I,\J)]^{-1} A(\I,j).
  \end{split}
 \end{equation}
 To compute the right--hand side we use only the elements of 
  selected columns $A(i,\J^{(t)}),$ 
  and rows $A(\I^{(s)},j).$ 
 Other elements of $A$ are not required to construct $\tilde A$ and we can avoid calculating them.
 Thus, evaluation and storage of $\tilde A$ requires $(mr+nr-r^2)$ matrix elements and is more cost--efficient that work with the whole matrix $A$ if $r\ll \min(m,n).$
 Due to the shape of the \emph{locus} of computed entries, shown on Fig.~\ref{fig:mat}, this decomposition is known as 
   \emph{skeleton}~\cite{gantmacher-1959}, 
   \emph{pseudo---skeleton} (if the exact inverse is replaced with, say, pseudo--inverse)~\cite{gt-psa-1995}, or
   \emph{cross}~\cite{tee-cross-2000}.
 
 \subsection{Notation for matrices and submatrices}
 Equation~\eqref{eq:mat} is understood element--wisely, i.e. holds for all possible values of free indices $i$ and $j.$ 
 According to the matrix multiplication rule, the summation is performed over the summation indices from the sets $\I$ and $\J,$ that are repeated in the formula, cf. Einstein's summation convention~\cite{Einstein-relativitat-1916}.
 Notation $A(\I,\J)$ refers to a submatrix on the intersection of rows $\I$ and columns $\J,$
  mimicking the intuitive syntax of programming languages like Fortran90, Matlab, R and Julia, where a vector of indices can be passed into an array to select a subsection of it, e.g. \texttt{A(1:2,1:3)} for a $2\times 3$ leading submatrix of $A.$
 We can also use index sets $\II=\{1,\ldots,m\}$ and $\JJ=\{1,\ldots,n\}$ to refer to full columns and rows.
 For instance, the approximant $\tilde A$ in~\eqref{eq:mat} is a product of three matrices:
  \begin{itemize}
   \item $m\times r$ matrix of columns $A(\II,\J)=[A(i,j)]_{i\in\II, j\in\J}$;
   \item inverse of the $r\times r$ submatrix at the intersection $A(\I,\J)=[A(i,j)]_{i\in\I, j\in\J}$;
   \item $r\times n$ matrix of rows $A(\I,\JJ)=[A(i,j)]_{i\in\I, j\in\JJ}$.
  \end{itemize}
 Embracing this notation, we will keep the same letter $A$ for all three factors of the cross interpolation.
 Compared to the $CGR$ notation~\cite{gt-psa-1995,gtz-psa-1997} or $CUR$ notation~\cite{mahoney-cur-2006}, our notation in~\eqref{eq:mat} highlights that factors of the cross decomposition are submatrices of the given matrix $A,$ which distinguishes it from SVD, QR and LU factorisations.

  \subsection{Maximum volume principle}
 The approximation $A\approx\tilde A$ is exact on the positions of computed rows $\I$ and columns $\J,$ which is why we call it \emph{interpolation}.
 For other entries the mismatch between $A$ and $\tilde A$ can be arbitrary large in general, because the approximation $\tilde A$ does not use any information about the most of $A$  apart of its few columns and rows.
 Theoretical error upper bounds can be obtained based on additional properties of the matrix, e.g. when $A=[f(x_i,y_j)]_{i,j=1}^{m,n}$ is generated by asymptotically smooth function~\cite{tee-mosaic-1996}.
 However, the quality of the cross approximation $\tilde A$ depends critically on a choice of good positions $(\I,\J)$ for the cross.
 Good theoretical estimates are available for the maximum--volume cross, i.e. such that $A(\I,\J)$ has the largest possible \emph{volume} 
 $$
 \vol A(\I,\J) = |\det A(\I,\J)|
 $$ 
 of all submatrices of this size.
 The maximum--volume principle for matrix approximation was first proposed in~\cite{gtz-maxvol-1997,gtz-psa-1997,gostz-maxvol-2010}, and 
   the estimates were later generalised to other norms~\cite{schneider-cross2d-2010,gt-skel-2011},
   and rectangular submatrices~\cite{zo-maxvol-2017,mo-rectmaxvol-2018}.
 
 Unfortunately, the search for a maximum--volume submatrix is NP--hard~\cite{bartholdi-1982} and cheaper alternative algorithms are required for practical calculations with large matrices.

 \subsection{Practical algorithms for matrix cross interpolation}
 When matrix $A$ is available in full, reliable algorithms for low--rank approximation are available, such as
  the famous singular value decomposition (SVD)~\cite{golubkahan-1965},
  and faster rank--revealing QR~\cite{gu-rrqr-1995} and LU~\cite{pan-rrlu-2000} algorithms.
 However, these approaches are unfeasible for very large--scale matrices, e.g. those coming from high--dimensional problems, when even $\O(mn)$ costs become prohibitive.

 To compute a sufficiently good cross with sublinear costs, the \emph{incomplete cross approximation}~\cite{tee-cross-2000} algorithm was proposed, that increases the volume of the intersection matrix by alternating updates of rows $\I$ and columns $\J.$
 In the set of columns $\J$ is fixed, there is a combinatorial number of possible row sets $\I$ to compare.
 To keep costs feasible, rows $\I$ are updated one--by--one with a greedy algorithm first suggested by Donald Knuth~\cite{knuth-1985}.
 Greedy updates of rows, shown in Alg.~\ref{alg:maxvol}, continue until the volume is large enough. 
 Then rows are fixed and columns are updated, and the algorithm alternates until a significantly large volume is obtained as desired.
 The details of this $\maxvol$ algorithm for matrix cross interpolation are given in~\cite{tee-cross-2000,gostz-maxvol-2010}.

\begin{algorithm}[t]
 \caption{One step of the practical row selection algorithm~\cite{knuth-1985}} \label{alg:maxvol}
 \begin{algorithmic}[1]
  \Require Sets $(\I,\J)$ of the interpolation~\eqref{eq:mat}
   \State 
    $
    B(\II,\I) \gets A(\II,\J)[A(\I,\J)]^{-1}.
    $
    \Comment{ $m \times n$ matrix with $B(\I,\I)=I$ }
   \State $(i^\new,i^\dagger) \gets \arg\max_{(i,j)\in \II\times\I}|B(i,j)|$
    \Comment{ $(i^\star,i^\dagger)\in\II\times\I$}
  \Ensure Updated row set $\I \gets \I \cup \{i^\new\} \setminus \{i^\dagger\}$ with $\vol A(\I,\J)\gets \vol A(\I,\J)|B(i^\new,i^\dagger)|.$
 \end{algorithmic}
\end{algorithm}

 A conceptually simpler \emph{adaptive cross approximation} (ACA) algorithm~\cite{bebe-2000} follows a greedy optimisation approach by increasing the interpolating sets by one columns and row at a time. 
 It can be seen as a Gaussian elimination with partial column pivoting~\cite{golub-2013}, which is computationally cheap but may result in exponential amplification $2^r$ of the error.
 A more conservative complete pivoting is believed to be numerically stable~\cite{wilkinson-1961}, but involves a search through all matrix elements, and thus is more expensive.
 A good alternative is the rook pivoting~\cite{poole-1992}, which searches for a pivot $(i^\new,j^\new)$ that is dominant in its own row and columns:
 \begin{equation}\label{eq:rook}
  |A(i^\new,j^\new)-\tilde A(i^\new,j^\new)| \geq |A(i,j)-\tilde A(i,j)|, \qquad 
      \text{for all $(i,j)$ such that $i=i^\new$ or $j=j^\new$}
 \end{equation}
 Rook pivoting avoids exponential deterioration of error~\cite{foster-1997,poole-2000} and has in practice the same asymptotical complexity as partial pivoting, so it seems to combine the best of both worlds.
 We use rook pivoting in combination with random pivoting, as shown in Alg.~\ref{alg:mat}.

\begin{algorithm}[t]
 \caption{One step of the matrix cross interpolation algorithm} \label{alg:mat}
 \begin{algorithmic}[1]
  \Require Sets $(\I,\J)$ of the interpolation~\eqref{eq:mat}
    \State Pick a random set of samples $\Rnd=\{(i,j)\}$ and choose the one with the largest error,
      \par
      $
      (i^\new,j^\new) \gets \arg\max_{(i,j)\in\Rnd} |A(i,j) - \tilde A(i,j)|
      $
    \Repeat\Comment{column and row partial pivoting updates}
     \State 
      $
      (i^\new,j^\new) \gets \arg\max_{i\in\II} |A(i,j^\new) - \tilde A(i,j^\new)|
      $
     \State 
      $
      (i^\new,j^\new) \gets \arg\max_{j\in\JJ} |A(i^\new,j) - \tilde A(i^\new,j)|
      $
    \Until{rook condition~\eqref{eq:rook} is met \textbf{or} computational budget is exhausted}
    \Ensure Expanded index sets $\I \gets \I \cup \{i^\new\}$, $\J \gets \J \cup \{j^\new\}$
 \end{algorithmic}
\end{algorithm}
 
 \begin{remark}[Numerical complexity] \label{rem:mat1}
  If $|\Rnd|=\O(m+n),$ a single rank--one update step evaluates $\O(m+n)$ matrix entries and performs $(m+n)r$ additional operations. 
  Thus $r$ steps of Algorithm~\ref{alg:mat} produce the rank--$r$ interpolation~\eqref{eq:mat} using $\O((m+n)r)$ matrix elements  plus $\O((m+n)r^2)$ additional operations.
 \end{remark}
 
 \begin{remark}[Accuracy]\label{rem:mat2}
  Algorithm~\ref{alg:mat} does not access all elements of the matrix and therefore is \emph{heuristic}, i.e. its accuracy can not be guaranteed in general.
 \end{remark} 
 
 Algorithm~\ref{alg:mat} is written in a very general way and many details are clearly improvable.
 For example, the choice of $\Rnd=\{(i,j)\}$ for initial sampling can be optimised to ensure $i\notin\I$ and $j\notin\J$ since the error $A-\tilde A$ is zero on the positions of the cross.
 A variety of other heuristic tricks were proposed, e.g.
   Mahoney et al~\cite{mahoney-cur-2006} suggest to estimate the column and row norms of $A$ and sample $(i,j)\in\Rnd$ with probabilities proportional to these norms.
 The focus of this paper is not the `best heuristic' for the matrix case, but the extension to high--dimensional problems.
 We refer the reader to~\cite{kumar-review-2017} for the review of matrix low-rank approximation algorithms.
 
 \section{Cross approximation and cross interpolation in higher dimensions}
 \subsection{Notation for tensors and multi--indices}
 We consider an array $A=[A(i_1,\ldots,i_d)]$ with $d$ indices $i_k$, $k=1,\ldots,d,$ which are also called \emph{dimensions} or \emph{modes}.
 Each index assumes values $i_k\in\II_k=\{1,\ldots,n_k\},$ where $n_k$ is called the \emph{mode size}.
 Such arrays are called \emph{tensors} in numerical linear algebra (NLA) community~\cite{golub-2013}, although we do not differentiate upper and lower indices, as it is customary for tensors in mathematical physics~\cite{Einstein-relativitat-1916}.
 The total storage required for $A$ grows exponentially with the dimension, prohibiting work with full $A$ for large $d.$
 Hence, tensor product representations are required for all practical calculations with tensors.

 At the heart of tensor product formats lies the idea of \emph{separation of indices}.
 Consider grouping indices $i_1,\ldots,i_k$ together and separating them from the group $i_{k+1},\ldots,i_d,$ 
   thus reshaping $n_1\times n_2\times \cdots \times n_d$ tensor $A$ into a $(n_1\cdots n_k)\times(n_{k+1}\cdots n_d)$ matrix
 $$
 A^{\{k\}}(i_{\leq k},i_{>k}) = A^{\{k\}}(i_1i_2\ldots i_{k};i_{k+1}\ldots i_d) = A(i_1,i_2,\ldots,i_d),
 $$
 called $k$--th matricization or unfolding of the tensor.
 As before, the equation is understood element--wisely for all possible values of all indices, i.e. $A^{\{k\}}$ differs from $A$ only by `shape'.
 Rows and columns of $A^{\{k\}}$ are enumerated by multi--indices
 $$
  i_{\leq k}=i_1i_2\ldots i_{k}\in\II_1\times\II_2\times\cdots\times\II_k, \qquad
  i_{>k}=i_{k+1}\ldots i_d\in\II_{k+1}\times\cdots\times\II_d.
 $$
 To separate row and column (multi)--indices, we apply matrix interpolation formula~\eqref{eq:mat} to $A^{\{k\}},$ yielding
 $$
  A^{\{k\}}(i_{\leq k},i_{>k})
  \approx \tilde A^{\{k\}}(i_{\leq k},i_{>k})
  = A^{\{k\}}(i_{\leq k},\I_{>k}) [A^{\{k\}}(\I_{\leq k},\I_{> k})]^{-1} A^{\{k\}}(\I_{\leq k},i_{> k}).
 $$
 Here $(\I_{\leq k},\I_{>k})$ indicate the positions of $r_k$ rows and columns of the interpolation cross in the unfolding $A^{\{k\}}.$ 

 \subsection{Tensor train format}
 The use of element--wise notation allows us to drop the superscript for the unfolding, because the dimensions of matrices and tensors are given by the range of the variables within.
 Hence, the equation above can be simplified as
 $$
 A(i_1,\ldots,i_d)
 \approx 
  A(i_1,\ldots,i_k,\I_{>k}) [A(\I_{\leq k},\I_{> k})]^{-1} A(\I_{\leq k},i_{k+1},\ldots,i_d),
 $$ 
 that emphasises separation of left and right groups of indices.
 By continuing the separation process, we arrive to the decomposition where all $i_k$'s are isolated:
 \begin{equation}\label{eq:tt}
  \begin{split}
  A(i_1,\ldots,i_d)
   & \approx \tilde A(i_1,\ldots,i_d)\\
   & = A(i_1,\I_{>1}) [A(\I_{\leq 1},\I_{>1})]^{-1} A(\I_{\leq 1},i_2,\I_{>2}) [A(\I_{\leq 2},\I_{>2})]^{-1}  \cdots  A(\I_{\leq d-1},i_d).
  \end{split}
 \end{equation}
 This formula is a direct generalisation of skeleton/cross interpolation~\eqref{eq:mat} to tensor case and is therefore called skeleton/cross tensor decomposition~\cite{ot-ttcross-2010}.
 It is a particular case of a more general tensor train (TT) decomposition~\cite{osel-tt-2011}, which appears if the factors of the TT decomposition are constructed from \emph{fibers} $A(\I_{\leq k-1},i_k,\I_{>k})$ of the given tensor.
 TT decomposition is itself a particular case of more general Hierarchical Tucker (HT) decomposition~\cite{hk-ht-2009,gras-hsvd-2010}.
 Cross approximation algorithms are available for HT format~\cite{lars-htcross-2013,bg-htcross-2015}, as well as for more specialised tensor formats, including Tucker~\cite{ost-tucker-2008} and canonical polyadic decomposition~\cite{sav-rr-2009}.
 \begin{remark}[Compression]
  The right--hand side of~\eqref{eq:tt} involves $\sum_{k=1}^d r_{k-1}n_kr_k - \sum_{k=1}^{d-1}r_k^2=\O(dnr^2)$ entries%
   \footnote{In all complexity estimates we assume $n_1\sim n_2\sim\cdots\sim n_d\sim n$ and  $r_1\sim r_2\sim\cdots\sim r_{d-1}\sim r.$} 
  of the tensor $A.$
 \end{remark}
 
 In general, tensor cross decomposition~\eqref{eq:tt} is not an interpolation formula.
 The following result from~\cite[Theorem 4]{sav-qott-2014} provides the sufficient condition for~\eqref{eq:tt} to be called tensor cross interpolation.
 \begin{theorem}[Interpolation, see~\cite{sav-qott-2014}]\label{thm:ii}
  If the crosses $(\I_{\leq k},\I_{>k})$ are \emph{nested}:
  \begin{equation}\label{eq:nest}
   \I_{\leq k+1} \subset \I_{\leq k} \times \II_{k+1},\qquad
   \I_{>k} \subset \II_{k+1} \times \I_{> k+1}, \qquad k=1,\ldots,d-1,
  \end{equation}
  formula~\eqref{eq:tt} interpolates the evaluated entries of the tensor,
  $$
  A(\I^{\leq k-1},i_k,\I^{>k}) = \tilde A(\I^{\leq k-1},i_k,\I^{>k}),\qquad k=1,\ldots,d.
  $$
 \end{theorem}
 Theorem~\ref{thm:ii} can not be reversed, i.e. nestedness of indices is not necessary for the interpolation, as shown by the following.
 
 \begin{theorem}[Exact recovery of the exact--rank tensor]\label{thm:exact}
   If $\rank A^{\{k\}}=r_k$ for all $k=1,\ldots,d-1,$ and 
   all submatrices $A(\I_{\leq k},\I_{>k})$ are non--singular, 
   the formula~\eqref{eq:tt} recovers the original tensor exactly,
   $
   A(i_1,\ldots,i_d) = \tilde A(i_1,\ldots,i_d).
   $
 \end{theorem}
 This theorem was first proven in~\cite{ot-ttcross-2010} with the additional requirement of nestedness. 
 %

 If $A^{\{k\}}$'s are only approximately low--rank, the good choice of crosses $(\I_{\leq k},\I_{>k})$ is important to ensure accurate approximation in~\eqref{eq:tt}.
 If all $(\I_{\leq k},\I_{>k})$ are maximum--volume submatrices in respective unfoldings $A^{\{k\}},$ the lower accuracy bounds are extended from matrices~\cite{gtz-maxvol-1997,gtz-psa-1997,gostz-maxvol-2010} to the tensor case~\cite[Theorem 1]{sav-qott-2014}.
 Inspired by the idea of maximal volume, we will now discuss practical algorithms for computation of sufficiently good crosses for the tensor cross interpolation.

 \subsection{Practical algorithms for tensor cross interpolation}
  In this section we provide a brief overview of tensor cross interpolation algorithms for TT format and compare them.
  \subsubsection{ALS maxvol algorithm~\cite{ot-ttcross-2010}}
  The algorithm in the pioneering paper~\cite{ot-ttcross-2010} is a direct generalisation of the matrix cross interpolation algorithm from~\cite{tee-cross-2000} to the tensor case.
 Starting from some selection of crosses $(\I_{\leq k},\I_{>k}),$ it updates them one--by--one using the maximum--volume principle.
  The left--to--right sequence of updates, called \emph{sweep}, is shown in Alg.~\ref{alg:als}.
  It is followed by a similar right--to--left sweep and the algorithm sweeps back and forth through the TT cores until convergence.
  This pattern of updates is often referred to as ALS, coming from \emph{alternating least squares} or \emph{alternating linear scheme}, although the abbreviation is often applied in broader sense.

  \begin{remark}[Nestedness in Alg.~\ref{alg:als}]
   The nestedness condition~\eqref{eq:nest} is not preserved during the sweep in Alg.~\ref{alg:als}. 
   Consider the moment when the left--to--right sweep reaches position $k$ in the train and replaces previous $\I_{\leq k}$ with the updated rows $\I_{\leq k}^\new.$
   The nestedness
   $
   \I_{\leq k}^\new \subset \I_{\leq k-1}^\new \times \II_k
   $ 
   is ensured by construction, so the nestedness of rows is maintained from the left side until the current active core.
   However 
   $
   \I_{\leq k+1} \not\subset \I_{\leq k}^\new \times \II_{k+1}
   $ 
   in general, because $\I_{\leq k+1}$ have not yet been updated and the nestedness of rows in the right part of the train is lost.

   The nestedness is recovered when the sweep reaches the end of the train, so the output $\tilde A$ of Alg.~\ref{alg:als} interpolates the given tensor $A.$
  \end{remark} 
  The main limitation of this algorithm is that it can not update the ranks $r_k$ of the interpolation, and therefore its success relies on two assumptions, both of which are not easy to ensure in practice:
  \begin{enumerate}
   \item the ranks $r_k$ of the interpolation $\tilde A$ are not underestimated to ensure that a good accuracy $|A-\tilde A|$ is achievable; and
   \item the ranks $r_k$ of the interpolation $\tilde A$ are not overestimated and non--singular submatrices $A(\I_{\leq k},\I_{>k})$ can be chosen at the initialisation step.
  \end{enumerate}

  \begin{algorithm}[t]
   \caption{Left--to--right sweep of the ALS maxvol cross approximation algorithm~\cite{ot-ttcross-2010}} \label{alg:als}
   \begin{algorithmic}[1]
    \Require Sets $(\I_{\leq k},\I_{>k})$ of the interpolation~\eqref{eq:tt}
    \For{$k=1,\ldots,d-1$}
     \State $\I_{\leq k}^\new \gets \maxvol [A(\I_{\leq k-1}^\new i_k,\I_{>k})]$ \Comment{choose $r_k$ rows in $r_{k-1}n_k \times r_k$ matrix}
    \EndFor
    \Ensure Updated index sets $\I_{\leq k} \gets \I_{\leq k}^\new,$ $k=1,\ldots,d-1$
   \end{algorithmic}
  \end{algorithm}

  \subsubsection{DMRG maxvol algorithm~\cite{so-dmrgi-2011proc}}
  To allow rank adaptation, we can consider a \emph{superblock} $A(\I_{\leq k-1}i_k,i_{k+1}\I_{>k+1})$ seen as $r_{k-1}n_k \times n_{k+1}r_{k+1}$ matrix.
  If we can compute the superblock in full, its low--rank decomposition can be computed by standard algorithms e.g. SVD~\cite{golubkahan-1965}.
  This allows us to adapt the rank $r_k$ in accordance with the desired accuracy and compute the good interpolation sets $(I_{\leq k},\I_{>k})$ from the factors of SVD decomposition, as shown in Alg.~\ref{alg:dmrg}.
  \begin{algorithm}[t]
   \caption{Left--to--right sweep of the DMRG maxvol cross approximation algorithm~\cite{so-dmrgi-2011proc}} \label{alg:dmrg}
   \begin{algorithmic}[1]
    \Require Sets $(\I_{\leq k},\I_{>k})$ of the interpolation~\eqref{eq:tt}, accuracy threshold $\eps$
    \For{$k=1,\ldots,d-1$}
     \State $B \gets A(\I_{\leq k-1}i_k,i_{k+1}\I_{>k+1})$ \Comment{compute superblock as $r_{k-1}n_k \times n_{k+1}r_{k+1}$ matrix}
     \State $USV^T\gets\textrm{svd}_\eps(B)$ \Comment{compute truncated SVD with accuracy $\eps$}
     \State $r_k\gets\rank(USV^T);$ 
            $\I_{\leq k}^\new \gets \maxvol U;$
            $\I_{>k}^\new \gets\maxvol V$
    \EndFor
    \Ensure Updated index sets $(\I_{\leq k},\I_{>k}) \gets (\I_{\leq k}^\new,\I_{>k}^\new),$ $k=1,\ldots,d-1$
   \end{algorithmic}
  \end{algorithm}

  \emph{Density matrix renormalization group} (DMRG)~\cite{white-dmrg-1992}  and related \emph{matrix product states} (MPS)~\cite{fannes-mps-1992,klumper-mps-1993} algorithms were developed in quantum physics community to find the ground state of a quantum spin system.
  The ranks of the ground state are not known in advance, which makes the rank adaptation crucial for the success of the method.
  Then the DMRG/MPS format was rediscovered in numerical linear algebra as the TT format~\cite{osel-tt-2011}, it was applied to a variety of problems including
  signal processing~\cite{dks-ttfft-2012,sav-rank1-2012},
  partial and fractional differential equations~\cite{DoOs-dmrg-solve-2011,DKhOs-parabolic1-2012,rst-volterra-2014},
  modelling of ionospheric plasma~\cite{dst-fb-2014} and 
  simulation of NMR~\cite{sdwk-nmr-2014}.
  Tailoring DMRG framework to compute interpolation and integration of high--dimensional functions is yet another example of extreme power and flexibility of algorithms, which can be understood, analysed and applied beyond the boundaries of the area where they were discovered.
 \begin{remark}[Nestedness in Alg.~\ref{alg:dmrg}]
   Similar to previous algorithm, Alg.~\ref{alg:als} does not preserve nestedness~\ref{eq:nest} during the sweep, but recovers it at the end of each sweep.
   Therefore, the output of Alg.~\ref{alg:dmrg} interpolates the initial tensor on all positions $(\I_{\leq k-1},i_k,\I_{>k}),$ $k=1,\ldots,d.$
 \end{remark}
  Unfortunately, Alg.~\ref{alg:dmrg} is moderately expensive --- it evaluates $\O(d n^2 r^2)$ points of the given tensor and interpolates only $\O(d n r^2)$ of them.

  \subsubsection{DMRG greedy algorithm~\cite{sav-qott-2014}}
  \begin{algorithm}[t]
   \caption{Left--to--right sweep of the DMRG greedy cross interpolation algorithm~\cite{sav-qott-2014}} \label{alg:cross}
   \begin{algorithmic}[1]
    \Require Sets $(\I_{\leq k},\I_{>k})$ of the interpolation~\eqref{eq:tt}
    \For{$k=1,\ldots,d-1$}
     \State Apply Alg.~\ref{alg:mat} to the superblock $A(\I_{\leq k-1}^\new i_k,i_{k+1}\I_{>k+1})$ 
     seen as $r_{k-1}n_k \times n_{k+1}r_{k+1}$ matrix.
     \par Find a new pivot $(i_{\leq k}^\new,i_{>k}^\new)$
     \State $\I_{\leq k}^\new  \gets \I_{\leq k}\cup\{i_{\leq k}^\new\};$
            $\I_{>k}^\new \gets \I_{>k}\cup\{i_{>k}^\new\}$
    \EndFor
    \Ensure Updated index sets $(\I_{\leq k},\I_{>k}) \gets (\I_{\leq k}^\new,\I_{>k}^\new),$ $k=1,\ldots,d-1$
   \end{algorithmic}
  \end{algorithm}

  \begin{figure}[t]
   \begin{center}
    \begin{tikzpicture}[x=5mm,y=-5mm,inner sep=.2em]
     \tikzstyle{box}=[thick,rounded corners=.3em,draw=black,fill=blue,fill opacity=.05];
     \tikzstyle{circ}=[black,fill=red,fill opacity=.15];
     \filldraw[box,dotted,fill opacity=0.0]   (0,0) rectangle (32,2);
     \filldraw[box]   (0,0) rectangle (4,2);
     \filldraw[box]  (10,0) rectangle (16,2);
     \filldraw[box]  (16,0) rectangle (22,2);
     \filldraw[box]  (28,0) rectangle (32,2);
     \filldraw[circ] ( 4,1) circle[x radius=2.0, y radius=1.0];
     \filldraw[circ] (10,1) circle[x radius=2.5, y radius=1.0];
     \filldraw[circ] (16,1) circle[x radius=2.0, y radius=1.0];
     \filldraw[circ] (22,1) circle[x radius=2.5, y radius=1.0];
     \filldraw[circ] (28,1) circle[x radius=2.5, y radius=1.0];
     \node at (1.5, 1) {$i_1$}; 
     \node at (6.7, 1) {$\ldots$}; 
     \node at (13.2,1) {$i_{k}$}; 
     \node at (18.7,1) {$i_{k+1}$};
     \node at (25,1) {$\ldots$};
     \node at (31,1)   {$i_d$};
     \node at (4,1)  {$\I_{> 1} \, \I_{\leq 1} $};
     \node at (10,1) {$\I_{> k-1} \, \I_{\leq k-1} $};
     \node at (16,1) {$\I_{> k} \, \I_{\leq k} $};
     \node at (22,1) {$\I_{> k+1} \, \I_{\leq k+1} $};
     \node at (28,1) {$\I_{> d-1} \, \I_{\leq d-1} $};
    \end{tikzpicture}
   \end{center}
   \begin{center}
    \begin{tikzpicture}[x=5mm,y=-5mm,inner sep=.2em]
     \tikzstyle{box}=[rounded corners=.3em,draw=black,fill=blue,fill opacity=.2];
     \filldraw[box,thick,fill=black!20] (0, 0) rectangle (14,10);
     \filldraw[box,fill=red]      (8, 0) rectangle (10,10); \node[red,below] at(9,1){$\I_{> k}$};
     \filldraw[box,fill=blue]     (0, 6) rectangle (14, 8); \node[blue,right]at(1,7){$\I_{\leq k}$};
     \node[left]  at(0,5){$\I_{\leq k-1} \times \II_k$};
     \node[above] at(7,0){$\II_{k+1}\times\I_{>k+1}$};
     \coordinate(new)at(5,3);
     \draw[very thick](0,3)--(14,3);
     \draw[very thick](5,0)--(5,10);
     \node[rectangle,rounded corners=.3em,fill=blue!20!white,fill opacity=.8,text opacity=1,left]at($(new)+(-.5,0)$){ $i_{\leq k}^{\new}$};
     \node[rectangle,rounded corners=.3em,fill=red!20!white,fill opacity=.8,text opacity=1,above]at($(new)+(0,-.7)$){$i_{> k}^{\new}$};
     \fill[red!70!black](new) circle(.3);
    \end{tikzpicture}
   \end{center}
   \caption{Cross interpolation algorithm~\cite{sav-qott-2014} searches for a new pivot $(i_{\leq k}^\new,i_{>k}^\new)$ in each superblock $A(\I_{\leq k-1}i_k,i_{k+1}\I_{>k+1})$ } \label{fig:alg}
  \end{figure}
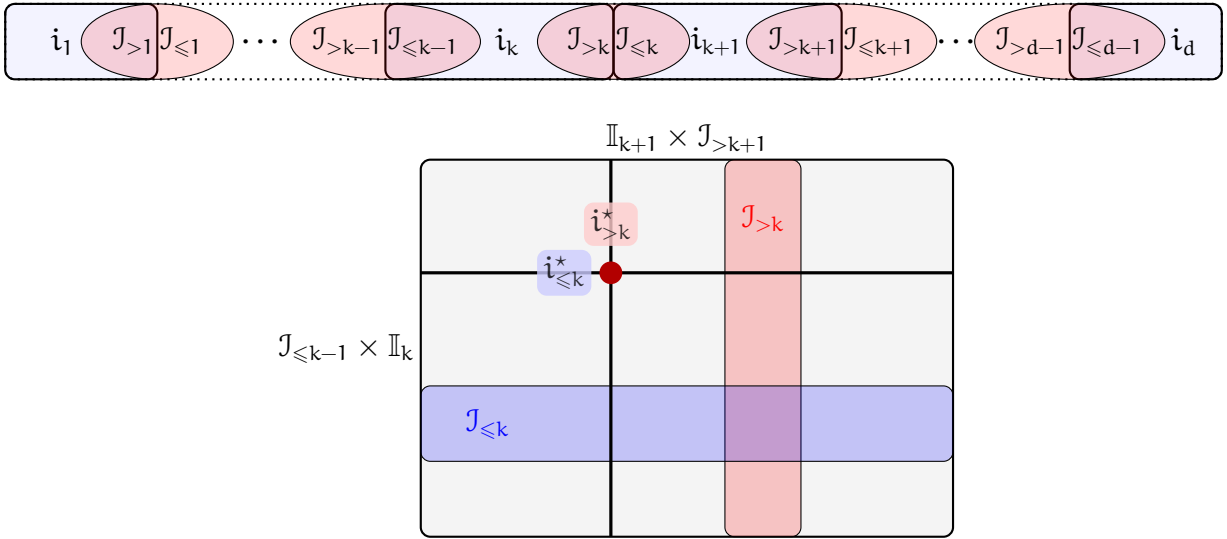

  Calculation of the superblock $A(\I_{\leq k-1}i_k,i_{k+1}\I_{>k+1})$ requires $\O(r^2n^2)$ function evaluations. 
  This may be too expensive, particularly when we aim for high precision and hence 
     employ large mode sizes $n_k$ for accurate quadratures and 
     expect large ranks $r_k$ to achieve accurate interpolation~\eqref{eq:tt}.
  To reduce costs we can replace $\maxvol$ optimisation step by greedy cross interpolation step, as proposed in~\cite{sav-qott-2014} and shown in Alg.~\ref{alg:cross} and Fig.~\ref{fig:alg}.
  The algorithm sweeps back and forth the tensor train~\eqref{eq:tt} and attempts to add one cross to each set $(\I_{\leq},\I_{>k})$ at a time.

 \begin{remark}[Nestedness in Alg.~\ref{alg:cross}]
   By construction, Alg.~\ref{alg:cross} preserves nestedness~\ref{eq:nest} at each internal step of the sweep.
   The output of Alg.~\ref{alg:cross} interpolates the initial tensor on all positions $(\I_{\leq k-1},i_k,\I_{>k}),$ $k=1,\ldots,d.$
 \end{remark}

 Alg.~\ref{alg:cross} requires $\O(dnr^2)$ evaluations of tensor elements and $\O(dnr^3)$ additional operations, which makes it one of the fastest tensor interpolation algorithms currently available in public domain.
 As all other algorithms considered in this section, it allows trivial parallelisation along each mode, which means that $\O(n)$ tensor entries forming each fiber can be evaluated in parallel.
 However, the fact that Alg.~\ref{alg:cross} maintains nestedness on each internal step makes it also suitable for parallelisation over all modes: 
 since no particular step can break the nestedness, all rank--one updates can be performed in parallel.
 This is explained in the following section.

 \subsection{Dimension parallel tensor cross interpolation algorithm}
  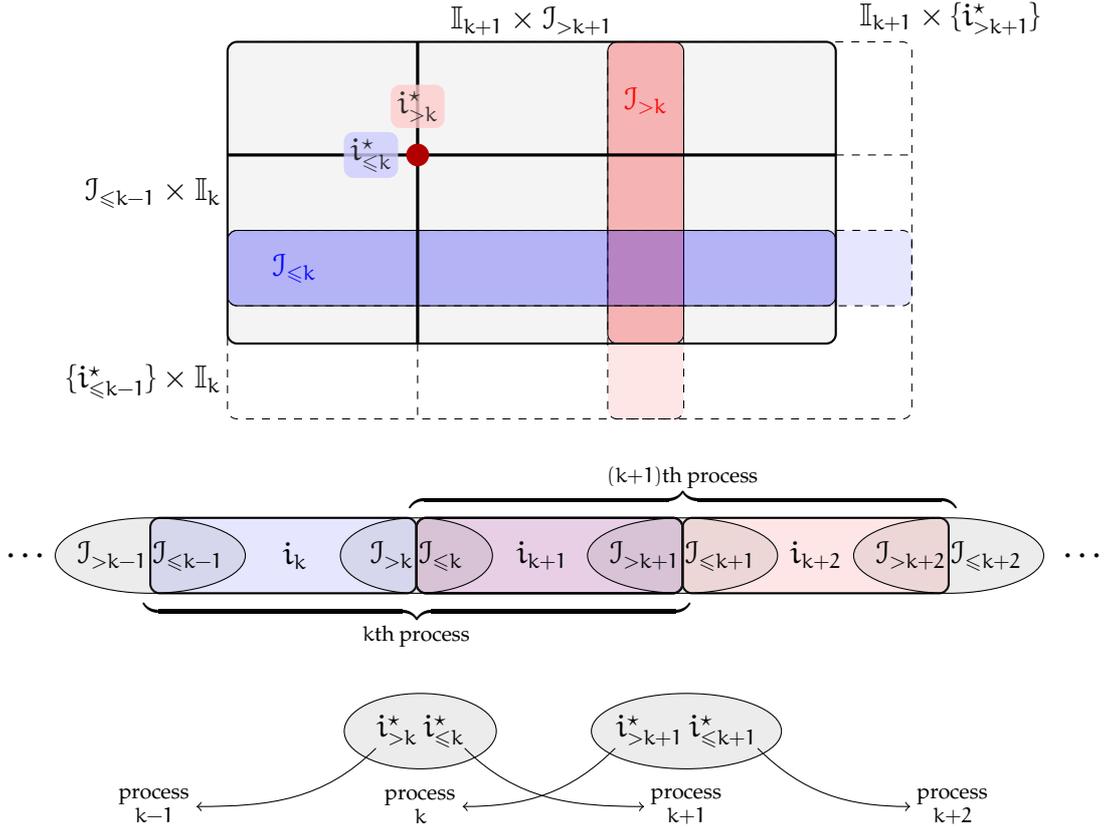
\begin{figure}[t]
   \begin{center}
    \begin{tikzpicture}[x=5mm,y=-5mm,inner sep=.2em]
     \tikzstyle{box}=[rounded corners=.3em,draw=black,fill=blue,fill opacity=.2];
     \tikzstyle{circ}=[black,fill=red,fill opacity=.10];
     \draw[box,dashed,fill opacity=0](0,0) rectangle (18,10);
     \draw[box,fill=red, dashed,fill opacity=.1](10,0) rectangle (12,10);
     \draw[box,fill=blue,dashed,fill opacity=.1](0,5) rectangle (18,7);
     \filldraw[box,thick,fill=black!20] ( 0, 0) rectangle (16,8);
     \filldraw[box,fill=red]      (10, 0) rectangle (12,8); \node[below,red] at(11,1){$\I_{> k}$};
     \filldraw[box,fill=blue]     ( 0, 5) rectangle (16, 7);\node[right,blue] at(1,6){$\I_{\leq k}$};
     \node[left]  at(0,4){$\I_{\leq k-1} \times \II_k$};
     \node[above] at(8,0){$\II_{k+1}\times\I_{>k+1}$};
     \node[left]  at(0,9){$\{i^{\new}_{\leq k-1}\}\times \II_k$};
     \node[above] at(19,0){$\II_{k+1}\times \{i^{\new}_{>k+1}\}$};

     \coordinate(new)at(5,3);
     \draw[very thick](0,3)--(16,3); \draw[dashed](16,3)--(18,3);
     \draw[very thick](5,0)--(5,8);  \draw[dashed](5,8)--(5,10);
     \node[rectangle,rounded corners=.3em,fill=blue!20!white,fill opacity=.8,text opacity=1,left]at($(new)+(-.5,0)$){ $i_{\leq k}^{\new}$};
     \node[rectangle,rounded corners=.3em,fill=red!20!white,fill opacity=.8,text opacity=1,above]at($(new)+(0,-.7)$){$i_{> k}^{\new}$};
     \fill[red!70!black](new) circle(.3);
    \end{tikzpicture}
   \end{center}
   \begin{center}
    \begin{tikzpicture}[x=5mm,y=-5mm,inner sep=.2em]
     \tikzstyle{boxb}=[thick,rounded corners=.3em,draw=black,fill=blue,fill opacity=.1];
     \tikzstyle{boxr}=[thick,rounded corners=.3em,draw=black,fill=red,fill opacity=.1];
     \tikzstyle{circ}=[black,fill=gray,fill opacity=.15];
     \filldraw[boxb]  (10,0) rectangle (17,2);
     \filldraw[boxb]  (17,0) rectangle (24,2);
     \filldraw[boxr]  (17,0) rectangle (24,2);
     \filldraw[boxr]  (24,0) rectangle (31,2);
     \filldraw[circ] (10,1) circle[x radius=2.5, y radius=1.0];
     \filldraw[circ] (17,1) circle[x radius=2.0, y radius=1.0];
     \filldraw[circ] (24,1) circle[x radius=2.5, y radius=1.0];
     \filldraw[circ] (31,1) circle[x radius=2.5, y radius=1.0];
     \node at (6.7, 1) {$\ldots$};
     \node at (13.8,1) {$i_{k}$};
     \node at (20.3,1) {$i_{k+1}$};
     \node at (27.5,1) {$i_{k+2}$};
     \node at (10,1) {$\I_{> k-1} \, \I_{\leq k-1} $};
     \node at (17,1) {$\I_{> k} \, \I_{\leq k} $};
     \node at (24,1) {$\I_{> k+1} \, \I_{\leq k+1} $};
     \node at (31,1) {$\I_{> k+2} \, \I_{\leq k+2}$};
     \node at (34.5,1) {$\ldots$};
     \node[anchor=north] at (17,2) {$\underbrace{\hspace*{17em}}_{k\text{th process}}$};
     \node[anchor=south] at (24,0) {$\overbrace{\hspace*{17em}}^{(k+1)\text{th process}}$};
    \end{tikzpicture}
   \end{center}
   \begin{center}
    \begin{tikzpicture}[x=5mm,y=-5mm,inner sep=.2em]
     \tikzstyle{circ}=[black,fill=gray,fill opacity=.15];
     \filldraw[circ] (17,1) circle[x radius=2.0, y radius=1.0];
     \filldraw[circ] (24,1) circle[x radius=2.5, y radius=1.0];
     \node[inner sep=0pt] (p1c) at (17,1) {$i_{> k}^{\new} \, i_{\leq k}^{\new} $};
     \node[inner sep=0pt] (p2c) at (24,1) {$i_{> k+1}^{\new} \, i_{\leq k+1}^{\new} $};
     \node (p1i) at (17,3) {$\substack{\text{process}\\ k}$};
     \node (p2i) at (24,3) {$\substack{\text{process}\\ k+1}$};
     \node (p0i) at (10,3) {$\substack{\text{process}\\ k-1}$};
     \node (p3i) at (31,3) {$\substack{\text{process}\\ k+2}$};
     \draw[->] (p1c.south east) to [out=-45,in=180] (p2i.west);
     \draw[->] (p2c.south west) to [out=-135,in=0] (p1i.east);
     \draw[->] (p2c.south east) to [out=-45,in=180] (p3i.west);
     \draw[->] (p1c.south west) to [out=-135,in=0] (p0i.east);
    \end{tikzpicture}
   \end{center}
   \caption{Parallel version of the cross interpolation algorithm. Top: excluding $i_{\leq k-1}^{\new}$ and $i_{>k+1}^{\new}$ from the row and column sets disentangles different steps in Alg. \ref{alg:cross}, while the new pivot search might be affected only a little or not at all if the pivot is located in the product of old subsets. Bottom: searching of pivots in different superblocks in parallel implies local data overlap and communication.} \label{fig:par}
  \end{figure}

Traditional ALS algorithm is carried out sequentially over tensor factors.
However, it was noticed that this dependence is more technical than essential.
A concurrency in ALS type algorithms is a matter of active research.
It was observed \cite{white-parallel-dmrg-2013} that the DMRG algorithm for ground state computations can be executed in parallel over subsets of TT blocks with only a little deterioration of the convergence.
Later a dimension parallel version of the HT-ALS for linear equations was developed \cite{etter-par-als-2016}.
In a non-adaptive HT Cross method the samples and the factors can also be reconstructed in parallel \cite{grasedyck-par-cross-2015}.

In this section we show that the adaptive Alg.~\ref{alg:cross} allows a natural parallelisation over dimensions.
From Line 3 of Alg.~\ref{alg:cross} we see that two consecutive steps $k$ and $k+1$ are connected by only one new pivot $i_{\leq k}^\new$, which expands the left index set $\I_{\leq k}$.
We can admit a slight restriction of the search space and replace expanded index sets $\I_{\leq k-1}^\new$ with the sets $\I_{\leq k-1}$ taken from the previous sweep, see Fig. \ref{fig:par} (top).
This restriction might potentially lead to a different (sub-optimal) pivot selection.
However, we observed no noticeable difference in the numerical experiments.
On the other hand, this allows us to search for new pivots in Line 2 of Alg. \ref{alg:cross} in a superblock $A(\I_{\leq k-1}i_k,i_{k+1}\I_{>k+1})$ with the old index sets, which is embarrassingly parallel over different $k$.
Different processes find their new pivots $(i_{\leq k}^\new,i_{>k}^\new)$ independently, communicate them and expand index sets before the next whole \emph{sweep} (instead of each next \emph{step} $k$ as in Alg. \ref{alg:cross}).
Since the superblocks owned by different processes overlap only for the neighbouring processes (e.g. the index $i_k$ belongs to only $(k-1)$th and $k$th superblocks),
only the neighbouring processes need to communicate:
the multi-index $i^{\new}_{\leq k}$ is sent from $k$th to $(k+1)$th process, and
$i^{\new}_{>k+1}$ is sent from $(k+1)$th to $k$th process, see Fig. \ref{fig:par} (bottom).

If fewer than $d-1$ processes are available,
each process can be given several consecutive superblocks.
The algorithm becomes similar to parallel DMRG \cite{white-parallel-dmrg-2013}, see Alg.~\ref{alg}:
each process performs the sequential sweep as in Alg. \ref{alg:cross} over its local chunk of the TT decomposition, and after that the neighbouring processes exchange new pivots in exactly the same way as described above.

\begin{remark}\label{rem:openmp}
This dimension parallel procedure can be hybridised with multi-threaded local computations, which consist of the evaluation of different samples in Alg. \ref{alg:mat} and the linear algebra of updating and applying the inversions $\left[A(\I_{\leq k}, \I_{>k})\right]^{-1}$.
\end{remark}

Assuming balanced splitting over $P$ processes, we conclude that each process performs $\O(dnr^2/P)$ evaluations of tensor elements and $\O(dnr^3/P)$ additional floating point operations.
Moreover, the tuples $i^{\new}_{\leq k}$, $i^{\new}_{>k}$ consist of at most $d-1$ integers, which need to be communicated with neighbours using $2$ messages in each of $r$ iterations, resulting in a total communication volume of $\O(dr)$.
Convergence checks require a global communication between all processors, amounting to $\O(r \log P)$ single--word messages in total.

  \begin{algorithm}[t]
   \caption{Dimension parallel DMRG greedy cross interpolation algorithm} \label{alg}
   \begin{algorithmic}[1]
    \Require Sets $(\I_{\leq k},\I_{>k})$ of the interpolation~\eqref{eq:tt}
    \State Deduce the range $[k_{\text{beg}}, k_{\text{end}})$ of superblocks belonging to the process $p$.
    \For{$k=k_{\text{beg}}, \ldots, k_{\text{end}} -1$} \Comment{in parallel over $p$}
     \State Apply Alg.~\ref{alg:mat} to the superblock $A(\I_{\leq k-1}^\new i_k,i_{k+1}\I_{>k+1})$
     seen as $r_{k-1}n_k \times n_{k+1}r_{k+1}$ matrix.
     \par Find a new pivot $(i_{\leq k}^\new,i_{>k}^\new)$
     \State $\I_{\leq k}^\new  \gets \I_{\leq k}\cup\{i_{\leq k}^\new\};$
            $\I_{>k}^\new \gets \I_{>k}\cup\{i_{>k}^\new\}$
    \EndFor
    \State Send $i_{\leq k_{\text{end}}-1}^\new$ to process $p+1$, receive $i_{\le k_{\text{beg}}-1}^\new$ from process $p-1$.
    \State Send $i_{>k_{\text{beg}}}^\new$ to process $p-1$, receive $i_{>k_{\text{end}}}^\new$ from process $p+1$.
    \State Update $\I_{\leq k_{\text{beg}}-1}^\new  \gets \I_{\leq k_{\text{beg}}-1}\cup\{i_{\leq k_{\text{beg}}-1}^\new\};$ $\I_{>k_{\text{end}}}^\new \gets \I_{>k_{\text{end}}}\cup\{i_{>k_{\text{end}}}^\new\}$.
    \Ensure Updated index sets $(\I_{\leq k},\I_{>k}) \gets (\I_{\leq k}^\new,\I_{>k}^\new),$ $k=1,\ldots,d-1$
   \end{algorithmic}
  \end{algorithm}

The parallelisation over the modes proposed in Alg.~\ref{alg} can scale well for the number of processes $P\lesssim d.$
It requires only a small number of global communications and lends itself well to distributed--memory `cluster' architectures and MPI--based implementation.
In contrast, the parallelisation along each mode requires all workers to access the shared block of memory where the fiber (or superblock) is stored.
Hence, this level of parallelisation is best for shared--memory architectures, such as cores and/or threads of a CPU/GPU processor and OpenMP--based implementation.
It scales efficiently when the number of cores/threads sharing the same memory is $T \lesssim n.$

In our algorithm we combine both of these approaches to achieve the best performance.

\section{High--dimensional integration}
 In this section we review quadrature rules for the numerical integration in high dimensions.
 We aim at computing an integral
 $$
 \Int = \int_{[0,1]^d} f(x_1,\ldots,x_d) \dx_1\cdots \dx_d = \int_{[0,1]^d} f(\x)\d\x,
 $$
 of a continuous function $f(\x)$ on a rectangular domain $[0,1]^d$.
 The exact integral is approximated by a quadrature
 $$
 \Int \approx \tilde\Int = \sum_{i=1}^{\neval} w_i f(\x_i),
 $$
 where $\neval$ nodes $\{\x_i\}$ and weights $\{w_i\}$ are properly chosen, such that the error $|\Int-\tilde\Int|$ is sufficiently small.
 Below we consider several examples of the quadrature rules.

 \subsection{Tensor product quadratures}\label{sec:quad}
 One of the simplest strategies is to rely on an appropriate one-dimensional quadrature rule (e.g. Gauss--Legendre, tahn-sin), 
   defined by the nodes $\{t_i\}_{i=1}^{n} \subset [0,1]$ and the weights $\{w_i\}_{i=1}^{n}$.
 The tensor product quadrature approximates each of the one-dimensional integrals independently,
 \begin{equation} \label{eq:quad_product}
  \tilde\Int = \sum_{i_1=1}^n \cdots \sum_{i_d=1}^n w_{i_1} \cdots w_{i_d} \cdot f(t_{i_1},\ldots,t_{i_d}).
 \end{equation}
 The main advantage of the tensor product quadrature is the fast convergence in $n$, which stems from the fast convergence of the one-dimensional Gauss--Legendre rule.
 For example, if a function $f(x)$, $x \in [-1,1]$,  is analytically extensible to a Bernstein ellipse $\mathcal{E}_{\rho}=\{z \in \mathbb{C}:~|z-1|+|z+1| \le \rho+\frac{1}{\rho}\}$ of radius $\rho>1$,
 the Gauss--Legendre quadrature converges with an exponential rate, $|\Int-\tilde\Int| = \O(\rho^{-n})$~\cite{Tadmor-exp_acc_diff-1986}.
 However, direct application of \eqref{eq:quad_product} is prohibitively expensive in high dimensions,
 as the total number of quadrature nodes $\neval=n^d$ grows exponentially with $d$.
 To utilise the benefits of the Gauss--Legendre quadrature in this case, we employ the TT approximation of $f(t_{i_1},\ldots,t_{i_d})$,
 which allows us to compute the quadrature with a linear cost with respect to $d$.
 Indeed, if we manage to separate function variables into the TT form~\eqref{eq:tt} as follows
 \begin{equation}\nonumber
  \begin{split}
   f(t_{i_1},\ldots,t_{i_d}) 
           & = A(i_1,\ldots,i_d) \\
           & \approx \tilde A(i_1,\ldots,i_d) 
             = A(i_1,\I_{>1}) [A(\I_{\leq 1},\I_{>1})]^{-1} A(\I_{\leq 1},i_2,\I_{>2}) \cdots  A(\I_{\leq d-1},i_d),
  \end{split}
 \end{equation}
 then plugging this in \eqref{eq:quad_product} we will rearrange the summation and treat each mode individually.
 The result is now given as a product of $(2d-1)$ matrices:
 $$
 \tilde\Int = \left(\sum_{i_1=1}^{n} w_{i_1} A(i_1,\I_{>1})\right) [A(\I_{\leq 1},\I_{>1})]^{-1} \left(\sum_{i_2=1}^{n} w_{i_2}A(\I_{\leq 1},i_2,\I_{>2})\right) \cdots  \left(\sum_{i_d=1}^{n} w_{i_d}A(\I_{\leq d-1},i_d)\right).
 $$

 If the TT ranks are bounded by $r$ that depends logarithmically on the accuracy, $r = \O(\log\eps^{-1})$,
 we obtain a poly-logarithmic overall complexity of the TT quadrature, $\O(d \log^3(\eps^{-1}))$.


 An alternative approach, which we find to converge faster in practice, is to incorporate quadrature weights together with the function values and apply the cross interpolation algorithm to their product, i.e. to
 $$
 B(i_1,\ldots,i_d) = w_{i_1} \cdots w_{i_d} \cdot f(t_{i_1},\ldots,t_{i_d}).
 $$
 This often leads to lower TT ranks/error compared to the approximation of $f(t_{i_1},\ldots,t_{i_d})$ if the function has a complicated structure near the boundaries.
 In this case, the boundary elements, multiplied by small cumulative products of the quadrature weights, become less influential to both the quadrature and the cross interpolation algorithm.

 \subsection{Monte Carlo and quasi Monte Carlo techniques} \label{sec:mc}
The Monte Carlo quadrature is a statistical method which is based on the central limit theorem.
It introduces random nodes $\{\x_i\}_{i=1}^{\neval}$ drawn from a uniform distribution on $[0,1]^d$,
and the integral is approximated by an average of the values of the function at these nodes and all weights equal,
\begin{equation}
 \tilde\Int = \frac{1}{\neval} \sum_{i=1}^{\neval} f(\x_i).
 \label{eq:mc}
\end{equation}
The integration error depends on the variance of $f(\x)$ (treated as a random field after randomisation of the coordinates $x$),
$|\Int - \tilde\Int|^2 \le \frac{\mathrm{var}(f)}{\neval}$.
Provided that the variance is independent of the dimension, so is the error.
However, the decay rate of $\neval^{-0.5}$ is often prohibitively slow, especially if a high accuracy is needed.

Quasi Monte Carlo (qMC)~\cite{nieder-qmc-1978,morokoff-qmc-1995} is another family of equal--weight quadrature rules (that is, $w_i = 1/\neval$ for all $i=1,\ldots,\neval$),
but the nodes are chosen semi-deterministically.
Firstly, one constructs a deterministic lattice rule, defined by a generating vector $\qmc=(q_1,\ldots,q_d).$
The lattice is optimised to minimise the worst--case error component by component~\cite{graham-QMC-2011,Kuo-QMC-2013}.
The quadrature nodes are then computed as shifted multiples of the generating vector modulo the interval $[0,1]$ in each variable,
\begin{equation} \label{eq:qmc}
 \x_i = \mathrm{frac}\left(\frac{i}{\neval}\qmc + \shift\right), \qquad i=1,\ldots,\neval.
\end{equation}
 Here $\shift=(s_1,\cdots,s_d)$ is a vector of random shifts, distributed uniformly on $[0,1]$, 
 and $\mathrm{frac}(x)$ denotes the fractional part of $x.$
Standard qMC rules provide a convergence rate $\O(\neval^{-\gamma})$, with $0.5\leq\gamma\leq1.$
Under certain assumptions on the function, the rate can be proven to be close to $1$, and the constant to be independent of $d$.
There exist higher order qMC rules \cite{Schwab-HOQMC-2014} which can achieve faster convergence, but at a price of more sophisticated lattice construction algorithms and stronger assumptions on the function.

The shifts $\shift$ make the quadrature \eqref{eq:qmc} unbiased, and they also allow to estimate the quadrature error.
We repeat qMC experiments using the same generating vector $\qmc$ but $S$ different shifts. 
Thus we obtain $S$ sets of nodes~\eqref{eq:qmc}, and use~\eqref{eq:mc} to calculate the estimators $\tilde\Int_j,$ $j=1,\ldots,S.$
Now the error can be estimated as the empirical standard deviation,
\begin{equation} \label{eq:std}
  \eps \approx \frac{1}{\langle \tilde\Int \rangle}\sqrt{\frac{1}{S-1}\sum_{j=1}^{S} \left(\tilde\Int_j - \langle \tilde\Int \rangle\right)^2}, \qquad  
  \langle \tilde\Int\rangle =  \frac{1}{S}\sum_{j=1}^{S} \tilde\Int_j.
\end{equation}
For the MC experiment we employ the same procedure by just sampling different $\neval$ points.

\section{Numerical experiments}
 \subsection{Ising integrals}
\begin{figure}[t]
     \begin{center}
      \begin{tikzpicture}[x={(20:7mm)}, y={(160:6mm)}, z={(90:8mm)}, baseline={(0,0,0)}]
      \pgfmathsetseed{12345678}
       \begin{scope}[canvas is yx plane at z=0] 
        \draw[gray] (.5,.5) grid[step=1] (6.5,6.5);
       \end{scope}
       \foreach \y in{1,...,6}{
        \foreach \x in{1,...,6}{
          \pgfmathrandominteger{\a}{1}{100}
          \ifnum\a>50
             \draw[thick,red,-latex]  (\x,\y,-.3)--(\x,\y,+.3);
             \shadedraw[ball color=red,draw opacity=.8,fill opacity=.8,ultra thin] (\x,\y,0) circle (.6mm);
          \else
             \draw[thick,blue,-latex] (\x,\y,+.3)--(\x,\y,-.3);
             \shadedraw[ball color=blue,draw opacity=.8,fill opacity=.8,ultra thin] (\x,\y,0) circle (.6mm);
          \fi
         }
        }
      \end{tikzpicture} 
      \hspace{10mm}
      \begin{tikzpicture}[x={(20:7mm)}, y={(160:6mm)}, z={(90:8mm)}, baseline={(0,0,0)}]
      \pgfmathsetseed{12345678}
       \begin{scope}[canvas is yx plane at z=0] 
        \draw[gray] (.5,.5) grid[step=1] (6.5,6.5);
       \end{scope}
       \foreach \y in{1,...,6}{
        \foreach \x in{1,...,6}{
         \pgfmathparse{\x+\y<5 || \x+\y>10 || (\x==5 && \y==2)  ? 0 : 1}
          \ifnum\pgfmathresult=0
             \draw[thick,blue,-latex] (\x,\y,+.3)--(\x,\y,-.3);
             \shadedraw[ball color=blue,draw opacity=.8,fill opacity=.8,ultra thin] (\x,\y,0) circle (.6mm);
          \else
             \draw[thick,red,-latex]  (\x,\y,-.3)--(\x,\y,+.3);
             \shadedraw[ball color=red,draw opacity=.8,fill opacity=.8,ultra thin] (\x,\y,0) circle (.6mm);
          \fi
        }
       }
       \draw[thick,black,-latex] (3.5,7.3,-1.8) -- ++(0,0,2) node[right]{$H$};
      \end{tikzpicture} 
     \end{center}
   \caption{Two--dimensional Ising model shown as a square lattice of interacting spins. 
   Normally, one would expect to observe individual spins in both states $\spin_{i,j}\in\{\spinup,\spindown\}$ with equal probability (as on the left panel).
   Spins also would align with the direction of external magnetic field (as shown on the right panel).
   Surprisingly, ferromagnetics will also exhibit collective large--distance behaviour (e.g. spontaneous magnetisation) at $H=0$ for sub--critical temperatures $T<T_c.$
   Theoretical explanation of this fact was first proposed by Lars Onsager in 1944.
   }
   \label{fig:2d}
  \end{figure}
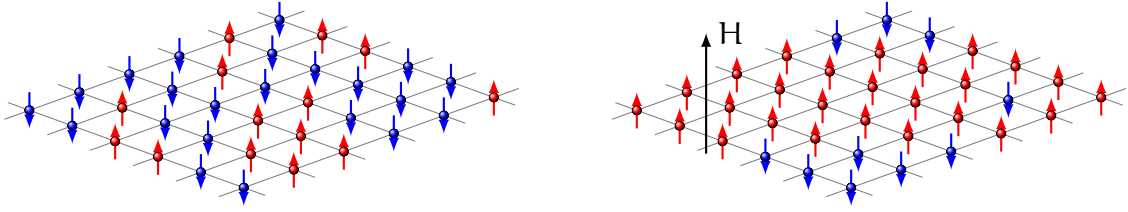

  To demonstrate the efficiency of the proposed approach, we apply tensor product interpolation to calculate high--dimensional integrals of so-called Ising class~\cite{bailey-ising-2006}.
  They are motivated by the famous 2D Ising model, explaining spontaneous magnetisation in ferromagnetic materials.
  It describes a ferromagnet as a rectangular $M\times N$ grid of spin--$\tfrac12$ particles where each spin $\spin_{i,j}$ can be observed in one of two possible states,
  $
   \spin_{i,j}\in\{+\tfrac12,-\tfrac12\} = \{\spinup,\spindown\}.
  $
  The energy of configuration $\spin=\{\spin_{i,j}\}_{\substack{i=1,\ldots,M\\ j=1,\ldots,N}}$ in magnetic field $H$ is given as follows:
  $$
   E(\spin) = \underbrace{-\sum_{i,j} \spin_{i,j}\spin_{i,j+1} 
                          -\sum_{i,j} \spin_{i,j}\spin_{i+1,j}}_{\text{next neighbour interaction}} 
            - \underbrace{H\sum_{i,j}\spin_{i,j}.}_{\mathclap{\text{response to magnetic field}}}
  $$
  The probability of each configuration is given by the Gibbs measure $\exp(-E(\spin)/kT)/Z,$
   where $T$ denotes the temperature and 
  $
  Z(T,H) = \sum_{\spin} \exp(-E(\spin)/kT)
  $ 
  is known as partition function.
  Assuming temperature and volume are constant, the Helmholtz free energy of the system is $F = -kT \log Z(T,H),$ and energy per particle is $f(T,H) = \lim_{\substack{M\to\infty \\ N\to\infty}} F(T,H)/(MN).$
  We may be interested in 
   \emph{spontaneous magnetisation} $m_0(T)=-\left.\frac{\d f}{\d H}\right|_{H=0}$
   and zero-field \emph{magnetic susceptibility} $\chi_0(T)=-\left.\frac{\d^2 f}{\d H^2}\right|_{H=0}.$
  Susceptibility is particularly interesting as it relates to long--distance spin--spin correlation and hence can explain collective behaviour in a ferromagnetic system which is connected by only next--neighbour interactions as shown in Fig.~\ref{fig:2d}.
  
  The 2D Ising model was first solved by Lars Onsager in 1944, who has never published the results.
  The solution for the magnetisation was published by Yang~\cite{yang-ising-1952}, and the susceptibility was calculated by Wu, McCoy, Tracy and Barouch~\cite{wmtb-ising-1976} as
  $$
  kT \chi_{0,\pm}(T) 
     = C_{0,\pm} \left| 1 - T/T_c \right|^{-7/4} 
     + C_{1,\pm} \left| 1 - T/T_c \right|^{-3/4}
     + \O(1),
  $$
  where $T_c$ denotes critical (Curie) temperature, which for the square and isotropic lattice is given by 
  $
  k T_c = 2/\ln(1+\sqrt2),
  $
  and $\pm$ refers to $T \to T_c$ from above ($+$) or below ($-$).
  The coefficients of the asymptotic expansion are given as infinite series,
  $$
        C_{0,+} \sim C_{1,+} \sim \sum_{\text{$d$ odd}}  \frac{\pi \D_d}{(2\pi)^d}, \qquad
        C_{0,-} \sim C_{1,-} \sim \sum_{\text{$d$ even}} \frac{\pi \D_d}{(2\pi)^d},
  $$
  where $\D_d$'s are $(d-1)$--dimensional integrals, which can be written as shown below~\cite{bailey-ising-2006}:
  \begin{align}
   \C_d & = 2 \int_{[0,1]^{d-1}} B_d(x_2,\ldots,x_d) \dx_2 \cdots \dx_d,                     \label{eq:c}   \\
   \D_d & = 2 \int_{[0,1]^{d-1}} A_d(x_2,\ldots,x_d) B_d(x_2,\ldots,x_d) \dx_2 \cdots \dx_d, \label{eq:d}   \\
   \E_d & = 2 \int_{[0,1]^{d-1}} A_d(x_2,\ldots,x_d) \dx_2 \cdots \dx_d,                     \label{eq:e}
  \end{align}
  with
  \begin{align*}
   A_d(x_2,\ldots,x_d) & = \prod_{1 \leq i < j \leq d} \left( \frac{1 - x_{i+1}\cdots x_j }{1 + x_{i+1} \cdots x_{j}} \right)^2, \\
   B_d(x_2,\ldots,x_d) & = \left(1+\sum_{k=2}^d x_2\cdots x_k \right)^{-1} \left(1+\sum_{k=2}^d x_k\cdots x_d \right)^{-1}.
  \end{align*}
  
  Bailey et al~\cite{bailey-ising-2006} took up a challenge to calculate $\D_d$'s numerically with high accuracy and then use inverse symbolic calculator~\cite{bailey-pslq-2000} to conjecture the values in closed form as a linear combination of physically relevant constants.
  Integrals $\C_d$ and $\E_d$ were motivated symbolically as a `simpler versions' of $\D_d$ in assumption that their values may also lead to certain insights.
  Indeed, all $\C_d$'s were analytically reduced to two--dimensional integrals and resolved numerically to extreme precision~\cite{bailey-ising-2006}. 
  Evaluation of $\D_d$'s and $\E_d$'s, even after significant analytic simplifications, proved to be difficult and accurate values were only obtained for relatively small dimensions. 
  We pick up the baton and consider the same problems, using the available values of $\C_d$ to verify the accuracy of the proposed tensor product algorithm, before proceeding to calculate $\D_d$'s with high accuracy for $d\lesssim1000.$

 \subsection{Experiment setup for double--, quadruple-- and high--precision calculations}
 Following Bailey~\cite{bailey-ising-2006}, we evaluate the integrals numerically using tensor product of one-dimensional Gauss--Legendre quadratures, as explained in Sec.~\ref{sec:quad}.
 The number of quadrature points in each direction, $n,$ is chosen adaptively to reach the desired accuracy.
 Since functions $A_d$ and $B_d$ are infinitely smooth, the Gauss--Legendre quadrature for $\C_d,$ $\D_d$ and $\E_d$ converges exponentially, and we can expect the number of accurate digits to grow linearly with $n.$

  The parallel implementation of the proposed algorithm is implemented in \textsc{Fortran} by authors.

  Double--precision calculations are implemented using  \textsc{GNU Fortran} compiler with \textsc{BLAS} and \textsc{Lapack} libraries from \textsc{Intel MKL}.

  For quadruple--precision calculations we compile the same code using a compiler option \texttt{-fdefault-real-8}, that sets the default size for \texttt{double precision} to $16$ bytes and increases precision to approximately $33$ decimal digits.
  We compiled the reference implementation of \textsc{BLAS} and \textsc{Lapack} libraries with the same parameter to reach quadruple precision in the whole calculation.

  For high--precision calculations we used the \textsc{MPFUN2015} library~\cite{bailey-highprec-2015,bailey-mpfun}.
  We had to rewrite reference implementation of necessary \textsc{BLAS} libraries to use the \texttt{mp\_real} data type offered by \textsc{MPFUN}.
  The code itself was compiled using the same compilers and options as for double precision calculations.
  The \textsc{MPFUN2015} library was set up to provide accuracy of $120$ decimal digits.

  The experiments were performed on two computers:
  \begin{itemize}
   \item at the University of Bath:      this research made use of the Balena High Performance Computing (HPC) Service.
   Each node on Balena contains an Intel Xeon E5-2650 v2 CPU with 16 cores, running at 2.6 GHz. A single job can occupy up to 32 nodes for 5 days.
  \item at the University of Brighton:  the development, testing and numerical experiments were made possible by use of a dedicated workstation.
   The workstation has two Intel Xeon E5-2650 v4 CPUs with 12 cores and 2 threads each, running at 2.2 GHz.
   It is also equipped with $0.5$ TB of operating memory, which proved essential for large--scale calculations reported below.
  \end{itemize}

 \subsection{Verification and benchmarking of the cross interpolation algorithm}
 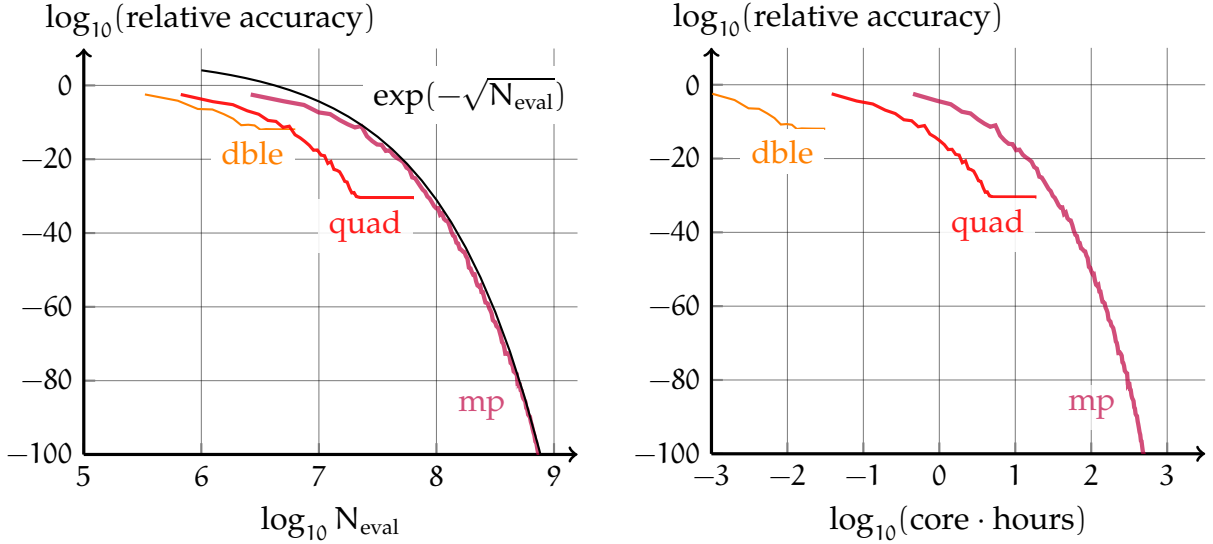
\begin{figure}[t]
   \hfil
  \begin{tikzpicture}
   \begin{axis}[
     xlabel={$\log_{10}\neval$},
     ylabel={$\log_{10}(\text{relative accuracy})$},
     xmin=5,xmax=9.2,
     ymax=10,ymin=-100,minor y tick num=0,
     ]
     \addplot[ttd,thick,opacity=1]        table[header=true,x=log10n,y=log10acc] {./ising/c1024_n33.dat} node[below left,lab]{dble};
     \addplot[ttq,very thick,opacity=.9]  table[header=true,x=log10n,y=log10acc] {./ising/c1024_n65.dat} node[below left,lab]{quad}; 
     \addplot[ttm,ultra thick,opacity=.7] table[header=true,x=log10n,y=log10acc] {./ising/c1024_n257_r200.dat} node[near end,below left,lab]{mp};
     \addplot[black,thick,domain=6:9] {8-0.0039*sqrt(10^x)} node [very near start,above right,lab]{$\exp(-\sqrt{\neval})$};
   \end{axis}
  \end{tikzpicture}
   \hfil
  \begin{tikzpicture}
   \begin{axis}[
     xlabel={$\log_{10}(\text{core}\cdot\text{hours})$},
     ylabel={$\log_{10}(\text{relative accuracy})$},
     xmin=-3,xmax=3.5,
     ymax=10,ymin=-100,minor y tick num=0,
     x filter/.code={\pgfmathparse{\pgfmathresult-log10(3600)}\pgfmathresult},
     ]
     \addplot[ttd,thick,opacity=1]        table[header=true,x=log10t,y=log10acc] {./ising/c1024_n33.dat} node[below left,lab]{dble};
     \addplot[ttq,very thick,opacity=.9]  table[header=true,x=log10t,y=log10acc] {./ising/c1024_n65.dat} node[below left,lab]{quad}; 
     \addplot[ttm,ultra thick,opacity=.7] table[header=true,x=log10t,y=log10acc] {./ising/c1024_n257_r200.dat} node[near end,below left,lab]{mp};
   \end{axis}
  \end{tikzpicture}
   \hfil

  \caption{Convergence of cross interpolation for calculation of $\C_{1024}$ in double, quadruple and multiple precision.
   Cross interpolation algorithm uses tensor product of one--dimensional Gaussian quadrature rules with $n=33$ points for double--precision, $n=65$ points for quadruple--precision and $n=257$ points for multiple--precision calculations.
   The results are verified against the $1000$--digit result reported in~\cite{bailey-ising-2006}. 
   The relative accuracy is shown w.r.t. number of function evaluations (left) and w.r.t. CPU time (right).
   We can clearly see that the proposed method converges exponentially.}
  \label{fig:mp}
 \end{figure}
 
 Bailey et al~\cite{bailey-ising-2006} found analytic transformation that converts $(d-1)$--dimensional integrals  $\C_d$ to two--dimensional form. 
 Using this two--dimensional representation, they calculated $\C_d$'s to $1000$ decimal digits for $d\leq1024.$
 They conjectured that
 $
 \C_{\infty} = \lim_{d\to\infty}\C_d = 2 e^{-2\gamma},
 $
 where $\gamma$ is the Euler--Mascheroni constant. 
 This result was later proven analytically.
 
 We compute $\C_{1024}$ directly as a $(d-1)$--dimensional integral using the proposed tensor product interpolation algorithm, and compare the numerical result with the one obtained by Bailey~\cite{bailey-ising-2006}.
 The comparison is shown at Fig.~\ref{fig:mp}.
 For double and quadruple precision calculations we observe an expected stagnation at the level of $15$ and $32$ decimal digits, respectively.
 When multiple precision calculations are used, the proposed algorithm seemingly provides exponential convergence for the integral $\C_{1024}.$
 As we can see on Fig.~\ref{fig:mp}, the observed convergence of relative accuracy $\eps$ agrees  well with the assumption
 $
 \eps \sim \exp(-\sqrt{\neval}). 
 $
 Since the number of samples evaluated by the cross interpolation algorithm is $\neval \sim dnr^2,$ 
  and $d,n$ remain constant, this allows us to conjecture that
 $
 \eps \sim \exp(-r),
 $ 
 i.e. the relative accuracy improves exponentially with the average TT-rank $r.$
 This makes tensor product decompositions preferable to currently known techniques such as MC and qMC algorithms.

 It should be noted that 
   although the use of quadruple and multiple precision calculations comes at a small extra cost in terms of number of points (it is sufficient to double the mode size $n$ to double the number of accurate digits), 
   it leads to significant overhead in terms of CPU time, since the quadruple and multiple precision calculations are not optimised to the same degree as native double precision calculations and \texttt{BLAS} libraries.
  This is why we report separately the convergence behaviour w.r.t. the number of evaluated points, and w.r.t. the CPU time on Fig.~\ref{fig:mp}.

 \subsection{Convergence and comparison with quasi Monte Carlo}
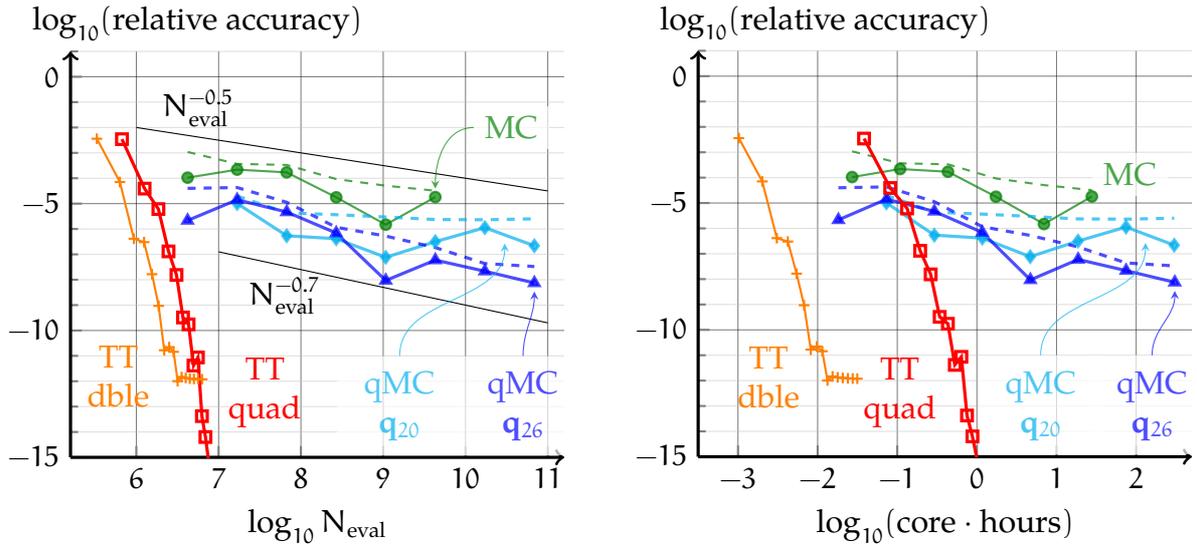
\begin{figure}[t]
   \hfil
  \begin{tikzpicture}
   \begin{axis}[%
     xlabel={$\log_{10}\neval$},
     ylabel={$\log_{10}(\text{relative accuracy})$},
     xmin=5.2,xmax=11.2,
     ymax=1,ymin=-15,minor y tick num=4,
     ]
     \addplot[mc,dashed]            table[header=true,x=log10n,y=log10sd]  {ising/mc/c1024.dat};
     \addplot[mc,mark=*]            table[header=true,x=log10n,y=log10acc] {ising/mc/c1024.dat} node(mc)[pos=1,outer sep=.5ex,inner sep=0ex]{};
     \addplot[qmc20,dashed]         table[header=true,x=log10n,y=log10sd]  {ising/qmc20/c1024.dat};
     \addplot[qmc20,mark=diamond*]  table[header=true,x=log10n,y=log10acc] {ising/qmc20/c1024.dat}node(qmc20)[pos=.9,outer sep=.5ex,inner sep=0ex]{};
     \addplot[qmc26,dashed]         table[header=true,x=log10n,y=log10sd]  {ising/qmc26/c1024.dat};
     \addplot[qmc26,mark=triangle*] table[header=true,x=log10n,y=log10acc] {ising/qmc26/c1024.dat} node(qmc26)[pos=1,outer sep=.5ex,inner sep=0ex]{};
     \addplot[ttd,mark=+]           table[header=true,x=log10n,y=log10acc] {ising/c1024_n33.dat} node[pos=0.95,left]{\begin{tabular}{c}TT \\ dble\end{tabular}};
     \addplot[ttq,mark=square]      table[header=true,x=log10n,y=log10acc] {ising/c1024_n65.dat} node[pos=0.35,right]{\begin{tabular}{c}TT\\ quad \end{tabular}};
     \addplot[no marks,black,domain=6:11] {1-0.5*x} node[pos=0.15,anchor=south]{$\neval^{-0.5}$};
     \addplot[no marks,black,domain=7:11] {-2-0.7*x} node[pos=0.2,anchor=north]{$\neval^{-0.7}$};
     
     \draw[stealth-,mc,thin]    (mc)    to[out=90,in=180] (axis cs:10.1,-2)node[right,lab]{MC};
     \draw[stealth-,qmc20,thin] (qmc20) to[out=270,in=90] (axis cs:9.2,-11)node[below,lab]{\begin{tabular}{c}qMC\\$\qmc_{20}$\end{tabular}};
     \draw[stealth-,qmc26,thin] (qmc26) to[out=270,in=90] (axis cs:10.7,-11)node[below,lab]{\begin{tabular}{c}qMC\\$\qmc_{26}$\end{tabular}};
   \end{axis}
  \end{tikzpicture}
   \hfil
  \begin{tikzpicture}
   \begin{axis}[%
     xlabel={$\log_{10}(\text{core}\cdot\text{hours})$},
     ylabel={$\log_{10}(\text{relative accuracy})$},
     xmin=-3.5,xmax=2.7,
     ymax=1,ymin=-15,minor y tick num=4,
     x filter/.code={\pgfmathparse{\pgfmathresult-log10(3600)}\pgfmathresult},
     ]
     \addplot[mc,dashed]            table[header=true,x=log10t,y=log10sd]  {ising/mc/c1024.dat};
     \addplot[mc,mark=*]            table[header=true,x=log10t,y=log10acc] {ising/mc/c1024.dat} node(mc)[pos=1,above right,lab]{MC};
     \addplot[qmc20,dashed]         table[header=true,x=log10t,y=log10sd]  {ising/qmc20/c1024.dat};
     \addplot[qmc20,mark=diamond*]  table[header=true,x=log10t,y=log10acc] {ising/qmc20/c1024.dat}node(qmc20)[pos=.9,outer sep=.5ex,inner sep=0ex]{};
     \addplot[qmc26,dashed]         table[header=true,x=log10t,y=log10sd]  {ising/qmc26/c1024.dat};
     \addplot[qmc26,mark=triangle*] table[header=true,x=log10t,y=log10acc] {ising/qmc26/c1024.dat} node(qmc26)[pos=1,outer sep=.5ex,inner sep=0ex]{};
     \addplot[ttd,mark=+]           table[header=true,x=log10t,y=log10acc] {ising/c1024_n33.dat} node[pos=0.95,left]{\begin{tabular}{c}TT \\ dble\end{tabular}};
     \addplot[ttq,mark=square]      table[header=true,x=log10t,y=log10acc] {ising/c1024_n65.dat} node[pos=0.35,left]{\begin{tabular}{c}TT\\ quad \end{tabular}};

     \draw[stealth-,qmc20,thin] (qmc20) to[out=270,in=90] (axis cs:0.8,-11)node[below,lab]{\begin{tabular}{c}qMC\\$\qmc_{20}$\end{tabular}};
     \draw[stealth-,qmc26,thin] (qmc26) to[out=270,in=90] (axis cs:2.2,-11)node[below,lab]{\begin{tabular}{c}qMC\\$\qmc_{26}$\end{tabular}};
   \end{axis}
  \end{tikzpicture}
   \hfil
   
  \caption{Integral $\C_{1024}$ calculated by TT cross interpolation (Alg.~\ref{alg}), Monte Carlo (MC), and quasi Monte Carlo (qMC). 
          Cross interpolation algorithm uses tensor product of one--dimensional Gaussian quadrature rules with $n=33$ points for double--precision and $n=65$ points for quadruple--precision calculations.
          QMC algorithm uses lattice generating vectors $\qmc_{20}$ and $\qmc_{26}$ minimising the worst--case error on $2^{20}$ and $2^{26}$ points, respectively.
          Solid lines: errors of numerical methods verified against the result of Bailey et al~\cite{bailey-ising-2006}.
          Dashed lines: relative standard deviation estimates~\eqref{eq:std} of MC and qMC with number of repetitions $S=16.$
  Left: relative accuracy w.r.t. different numbers of function evaluations~$\neval.$ 
  Right: relative accuracy w.r.t. total CPU time.
  }
  \label{fig:qmc}
 \end{figure}

  On Fig.~\ref{fig:qmc}  the proposed algorithm is compared with state of the art Monte Carlo (MC) and Quasi MC approaches (see Sec. \ref{sec:mc}).
  For the MC quadrature we use uniformly distributed samples on $[0,1]^d$.

  For the qMC algorithm a particular care must be taken when choosing the correct lattice.
  Frances Kuo's website\footnote{\url{http://web.maths.unsw.edu.au/~fkuo/}} provides a large collection of pre-generated lattices which were generated by optimising the worst case error with product weight parameters $\gamma_k = k^{-2},$ motivated by stochastic PDEs.
  For the integrals considered in this paper all variables seem to play similar role and we would prefer a lattice with equal weights.
  Hence we used the component by component algorithm from Dirk Nuyens's website\footnote{\url{https://people.cs.kuleuven.be/~dirk.nuyens/qmc-generators/}} 
  and constructed generating vectors $\qmc_{20}$ and $\qmc_{26}$ by minimising the worst case error on $2^{20}$ and $2^{26}$ points respectively.
  Notice that the lattice generated from $\qmc_{20}$ starts repeating when the number of points exceeds $2^{20},$ leading to a visible stagnation of the $\qmc_{20}$ quadrature error in Fig.~\ref{fig:qmc}.
  This is why we created lattice $\qmc_{26}$ which remains convergent and allows to scale the computations up to billions of points.
  It has to be noted that optimising a lattice is rather expensive --- the CBC algorithm took several days to produce $\qmc_{26}$ (this cost is not included in further analysis).

  As in the previous subsection, we calculate $\C_{1024}$ and compare our results against the 1000-digit accurate value computed in~\cite{bailey-ising-2006}.
  These errors are plotted on solid lines in Fig. \ref{fig:qmc}.
  We also show by dashed lines the relative empirical standard deviation for MC and qMC algorithms as described in \eqref{eq:std}.
  Notice that the true error exhibits a higher fluctuation for different $\neval$, although the overall convergence trend coincides with that for the standard deviation.


  We see that the MC method converges with the rate $\neval^{-0.5}$ as expected from the CLT,
  while the qMC method (with $\qmc_{26}$) exhibits a higher rate $\neval^{-0.7}$.
  The TT decomposition has a much richer approximation capacity, and provides a sub-exponential convergence, as shown also in Fig.~\ref{fig:mp}.
  When all calculations are performed in double precision, TT cross interpolation is always faster than MC and qMC methods.
  Switching to quadruple precision increases the TT time significantly, since we lose optimisations of Intel MKL, but the rapid convergence still makes it the fastest method for high accuracy.

 \subsection{Evaluation of Ising susceptibility integrals}

 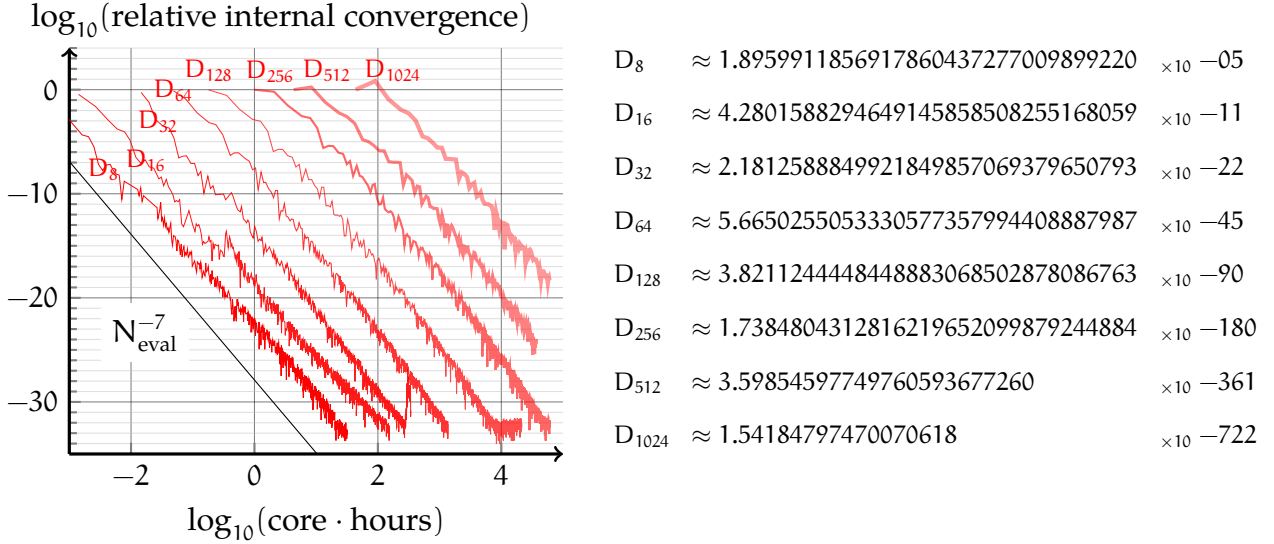
\begin{figure}[t]
  \begin{tikzpicture}
    \begin{axis}[
      xlabel={$\log_{10}(\text{core}\cdot\text{hours})$},
      ylabel={$\log_{10}(\text{relative internal convergence})$},
      xmax=5,xmin=-3,
      ymax=4,ymin=-35,minor y tick num=9,
      x filter/.code={\pgfmathparse{\pgfmathresult-log10(3600)}\pgfmathresult},
      ]
      \addplot[tt,ultra thin,opacity=1]  table[header=true,x=log10t,y=log10acc] {./ising/d8_n129.dat}   coordinate[pos=.02](8);
      \addplot[tt,very thin,opacity=.9]  table[header=true,x=log10t,y=log10acc] {./ising/d16_n129.dat}  coordinate[pos=.02](16);
      \addplot[tt,very thin,opacity=.9]  table[header=true,x=log10t,y=log10acc] {./ising/d32_n129.dat}  coordinate[pos=.01](32);
      \addplot[tt,thin,opacity=.8]       table[header=true,x=log10t,y=log10acc] {./ising/d64_n129.dat}  coordinate[pos=.00](64);
      \addplot[tt,opacity=.7]            table[header=true,x=log10t,y=log10acc] {./ising/d128_n129.dat} coordinate[pos=.00](128);
      \addplot[tt,thick,opacity=.6]      table[header=true,x=log10t,y=log10acc] {./ising/d256_n129.dat} coordinate[pos=.001](256);
      \addplot[tt,very thick,opacity=.5] table[header=true,x=log10t,y=log10acc] {./ising/d512_n65.dat}  coordinate[pos=.00](512);
      \addplot[tt,ultra thick,opacity=.4]table[header=true,x=log10t,y=log10acc] {./ising/d1024_n65.dat} coordinate[pos=.00](1024);
      
      \addplot[no marks,black,domain=0:5] {-3-7*x} node[midway,below left,lab]{$\neval^{-7}$};
      \node(d) at(axis cs:5.5,5){};
    \end{axis}
    \node[scale=.8,tt]                   at(8){$\D_{8}$};
    \node[scale=.8,tt]                   at(16){$\D_{16}$};
    \node[scale=.8,tt]                   at(32){$\D_{32}$};
    \node[scale=.8,tt]                   at(64){$\D_{64}$};
    \node[scale=.8,tt,above]             at(128){$\D_{128}$};
    \node[scale=.8,tt,above]             at(256){$\D_{256}$};
    \node[scale=.8,tt,above right]       at(512){$\D_{512}$};
    \node[scale=.8,tt,above right]       at(1024){$\D_{1024}$};
    \node[scale=.8,below right]at(d){
     $\begin{array}{lll}
      \D_{8}    &\approx  1.8959911856917860437277009899220&\e{-05}   \\[2ex]
      \D_{16}   &\approx  4.2801588294649145858508255168059&\e{-11}   \\[2ex]
      \D_{32}   &\approx  2.1812588849921849857069379650793&\e{-22}   \\[2ex]
      \D_{64}   &\approx  5.6650255053330577357994408887987&\e{-45}   \\[2ex]
      \D_{128}  &\approx  3.8211244448448883068502878086763&\e{-90}   \\[2ex]
      \D_{256}  &\approx  1.7384804312816219652099879244884&\e{-180}  \\[2ex]
      \D_{512}  &\approx  3.59854597749760593677260&\e{-361}  \\[2ex]
      \D_{1024} &\approx  1.54184797470070618&\e{-722} 
     \end{array}$
   };
   \end{tikzpicture}
  \caption{Evaluation of the Ising susceptibility integrals $\D_d$ given by~\eqref{eq:d}.
           The results are computed by the cross interpolation algorithm in quadruple precision using tensor product of one--dimensional Gaussian quadrature rules with $n=129$ points (for $\D_8$ to $\D_{256}$) and $n=65$ points (for $\D_{512}$ and $\D_{1024}$).
           Left: convergence of cross interpolation algorithm measured by the relative internal convergence, as a function of total CPU time spent on the calculation. 
           Right: values of the $\D_d$'s calculated by the proposed algorithm.}
  \label{fig:d}
 \end{figure}
 
 Now we attempt to compute original Ising susceptibility integrals $\D_d$ given by~\eqref{eq:d}.
 Computing $\D_d$'s for large $d$ is much more challenging than evaluating $\C_d$'s, for two reasons.
 Firstly, each evaluation of the integrand takes $\O(d)$ operations for $\C_d,$ but $\O(d^2)$ for $\D_d.$
 Secondly, all $\C_d$'s can be analytically reduced to two dimensional integrals, while for $\D_d$'s reduction performed in~\cite{bailey-ising-2006} only reduces the dimensionality by one in special cases.
 Using a combination of analytic transforms and Gaussian tensor--product quadratures, Bailey and collaborators calculated $\D_5$ to $500$ decimal digits using $18$h on 256 CPUs of IBM Power5 nodes at the Lawrence Berkeley National Laboratory.
 They also produced $\D_6$ to almost $100$ decimal digits.
 Using qMC algorithm, they also calculated $\D_7$ and $\D_8$ to $5$ decimal digits.
 Further integrals $\D_d$ were not made available.

 We apply the proposed tensor interpolation algorithm to calculate $\D_d$'s in the original form~\eqref{eq:d} as $(d-1)$--dimensional integrals.
 We use the quadruple--precision version of the code and aim to calculate integrals $\D_8,$ $\D_{16},$ $\D_{32},\ldots,\D_{1024}$ to about $30$ decimal digits, which is measured by the internal convergence.
 The convergence plots are shown on Fig.~\ref{fig:d}.
 The convergence rate is approximately of order $7$ for all considered integrals; noting a slight bent of the curve for $\D_{256}$ we are hopeful that exponential convergence could have been revealed if calculations were allowed to run longer and reach higher accuracy. 
  
 By looking at the values of $\D_d$'s on Fig.~\ref{fig:d} it is easy to note that they decay exponentially.
 This was noted by Bailey et al, who proved~\cite[Thm. 3]{bailey-ising-2006} that $\O(14^{-d})\leq\D_d\leq\O(4^{-d}).$
 They conjectured that as $d\to\infty,$ $\D_d\sim\Delta^{-d},$ and based on a few available for them values $\D_d$ estimated $\Delta\approx5.$
 Based on our values $\D_{128}$ and $\D_{256}$ shown in Fig.~\ref{fig:d}, we improve this estimate to
 \begin{equation}\label{eq:Delta}
  \Delta \approx 5.0792202086636783360436879567820.
 \end{equation}

  \subsection{Performance and scalability}
 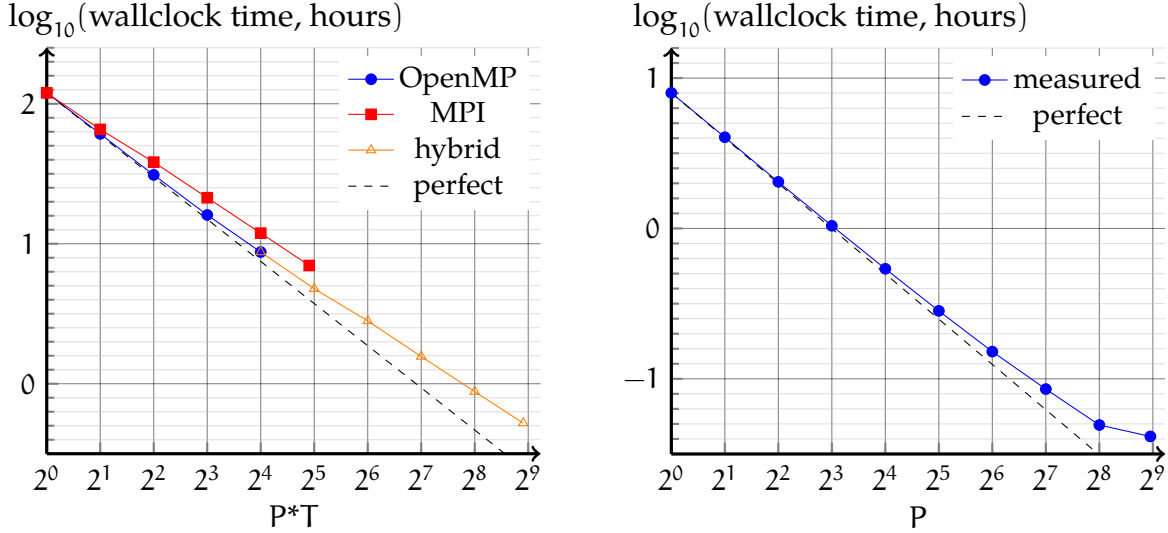
\begin{figure}[t]
   \hfil
   \begin{tikzpicture}
    \begin{axis}[
      xlabel={$P$*$T$},
      ylabel={$\log_{10}(\text{wallclock time, hours})$},
      xmode=log,
      xmin=1,xmax=600,
      xtick={1,2,4,8,16,32,64,128,256,512},
      xticklabels={$2^0$, $2^1$, $2^2$, $2^3$, $2^4$, $2^5$, $2^6$, $2^7$, $2^8$, $2^9$},
      ymax=2.4,ymin=-.5,ytick distance=1,minor y tick num=9,
      y filter/.code={\pgfmathparse{\pgfmathresult-log10(3600)}\pgfmathresult},
      ]
      \addplot[mark=*,blue]               table[header=true,x=np,y=log10t] {./ising/d32_openmp_scaling.dat}; \addlegendentry{OpenMP};
      \addplot[mark=square*,red]          table[header=true,x=np,y=log10t] {./ising/d32_mpi_scaling.dat}; \addlegendentry{MPI}; 
      \addplot[mark=triangle,orange]      table[header=true,x=np,y=log10t] {./ising/d32_hybrid_scaling.dat}; \addlegendentry{hybrid}; 
      \addplot[no marks,black,dashed,domain=1:512] {5.634-log10(x)}; \addlegendentry{perfect};
    \end{axis}
   \end{tikzpicture}
   \hfil
   \begin{tikzpicture}
    \begin{axis}[
      xlabel={$P$},
      ylabel={$\log_{10}(\text{wallclock time, hours})$},
      xmode=log,
      xmin=1,xmax=600,
      xtick={1,2,4,8,16,32,64,128,256,512},
      xticklabels={$2^0$, $2^1$, $2^2$, $2^3$, $2^4$, $2^5$, $2^6$, $2^7$, $2^8$, $2^9$},
      ymax=1.2,ymin=-1.5,ytick distance=1,minor y tick num=9,
      y filter/.code={\pgfmathparse{\pgfmathresult-log10(3600)}\pgfmathresult},
      ]
      \addplot[mark=*,blue]  table[header=true,x=np,y=log10t] {./ising/d512_mpi_scaling.dat}; \addlegendentry{measured}; 
      \addplot[no marks,black,dashed,domain=1:512] {4.458-log10(x)}; \addlegendentry{perfect};
    \end{axis}
   \end{tikzpicture}
   \hfil



  \caption{Left: strong scaling for $\D_{32}$ for different numbers of processes $P$ and numbers of threads $T$, quadruple precision with $n=129$. 
           Right: strong MPI scaling (1 OpenMP thread per process) for $\D_{512}$, double precision with $n=33$.
           }
  \label{fig:mpi16}
 \end{figure}
  In Fig.~\ref{fig:mpi16} we benchmark the algorithm for different numbers of processes and threads using MPI, OpenMP and hybrid parallelisation.
  The first two lines in Fig.~\ref{fig:mpi16} (left) show the CPU time for OpenMP-only parallelisation of local computations (i.e. essentially Alg.~\ref{alg:cross} with no dimension parallelisation), and for MPI-only approach where all local computations are performed in one thread,
  but different chunks of the TT decomposition are assigned to different processes (Alg.~\ref{alg}).
  Moreover, the hybrid approach always uses $T=16$ threads for local operations, and different numbers of processes $P$ for parallelisation over dimension.
  In Fig.~\ref{fig:mpi16} (left) we report the product of the number of processes and the number of threads in each process.

  Since the $\D_{32}$ integral involves actually a $31$-dimensional function, the maximal number of processes is limited by $30$.
  Here the hybrid framework allows us to accelerate the computing further up to a maximum of $512$ cores, available on the Balena cluster per one job.
  We notice a very good scaling, since the cost of communicating $\O(rd + r\log P)$ bytes is much smaller than the cost of computing $\O(dnr^2/P)$ tensor elements.
  A slight deviation from the linear scaling for the largest numbers of processes is due to load imbalance, as different TT blocks pick up different ranks in the course of the cross algorithm.

  This is demonstrated further in Fig.~\ref{fig:mpi16} (right), where we approximate a function for the $\D_{512}$ integral.
  The maximal number of processes $510$ allows us to use only $T=1$ OpenMP thread, and instead vary the number of MPI processes $P$ in the entire range.
  We see that the time is closer to the perfect scaling due to better balancing when each process owns more TT blocks.
  Even better scaling could be expected for $\D_{2^p+2}$ integrals, where the same number of TT blocks could be assigned to each of $2^p$ processes.
  Nevertheless, even in a deliberately unbalanced situation (which is more practical though),
  the algorithm scales almost linearly up to the maximum computing capacity available at the given machine.

 Finally, we should note that even though with the proposed algorithm~\ref{alg} we enjoy fast convergence, the numerical costs remain quite high.
 For example, calculation of $\D_{1024}$ to 18 decimal digits (see Fig.~\ref{fig:d}) took about 4 days on 512 nodes of Balena supercomputer at the University of Bath, consuming approximately a megawatt hour of energy.
 Based on our preliminary experiments with qMC, and assuming that the convergence rate $\eps\sim \neval^{-0.7}$ will not deteriorate, we estimate that to reach the same accuracy with qMC we would need approximately $10^{13}$ years of calculations and $10^9$ terawatt hours of energy ---
 which exceeds the age of the Universe ($\approx 1.3 \cdot 10^{10}$ years) and annual world energy consumption ($\approx 1.5\cdot 10^5$ Twh in 2014) by three orders of magnitude.

\section{Conclusion}
 The problem of high--dimensional integration is a particularly important and challenging area.
 Motivated by risk simulation in finance and engineering, this problem was actively researched 
   and resulted in Monte Carlo Metropolis algorithm~\cite{metropolis-mc-1949}, which is considered as one of top 10 algorithms of the \nth{20} century~\cite{dongarra-top10}.
 The use of random samples in the MC algorithm allows to break away from tensor--product quadratures and hence avoid the curse of dimensionality, seemingly inevitable in higher dimensions.
 The flexibility and simplicity of MC was spoiled by its slow convergence, motivating the further development, until the arrival of quasi Monte Carlo algorithms~\cite{nieder-qmc-1978,morokoff-qmc-1995}. 
 QMC methods can be optimised for a class of functions (e.g. those appearing from stochastic PDEs~\cite{graham-QMC-2011,Schwab-HOQMC-2014}), and demonstrate faster convergence, which currently makes them methods of choice in areas of sPDEs, finance and risk modelling, engineering, etc.
 However, the convergence is still not too fast, particularly considering that in practice many end users can make sub-optimal choices in choosing/creating the correct qMC lattice for their problems.
 
 The curse of dimensionality turns therefore in a challenge of precision. 
 Although admittedly many practical problems (e.g. in areas of stochastic inference or machine learning) do not require precision above one or two decimal digits, 
   many applications (e.g. engineering, theoretical quantum physics, quantum computations) need the answer to be precise to ten(s) or hundred(s) of decimal digits, which can't be achieved (or leads to excessive costs in terms of energy and CPU time) using mainstream MC/qMC approaches.
 In this paper we address this challenge by development of a new algorithm, based on tensor decompositions.
 We are pleased to see that the idea of the decompositional approach to matrix computation~\cite{stewart-lr-2000}, which was also recognised as a top 10 algorithm of \nth{20} century~\cite{dongarra-top10}, can break the curse of dimensionality --- arguably one of the main challenges of numerical mathematics since 1960s~\cite{bellman-dyn-program-1957} and till this day.
 Tensor product algorithms are undergoing very rapid development during the last $15$ years, both in terms of theory, algorithmic implementations, and applications.
 Using the idea of separation of variables, tensor methods give a new hope in lifting the curse of dimensionality and drastically reducing the computational burden associated with high--dimensional problems in a number of areas from quantum physics and chemistry to stochastics, signal processing and data analysis.
 In this paper we applied tensor cross interpolation algorithm~\cite{sav-qott-2014} to reconstruct the behaviour of the given high--dimensional function from a few samples and to numerically integrate it.
 Our research proposes a new step in development of tensor product algorithms, by combining 
   the algorithmic power provided by data--sparse low--rank tensor product representations, and 
   the efficient parallel implementation utilising the potential of modern HPC systems.
 
 The Ising susceptibility integrals, which we use in this paper to demonstrate the efficiency of the proposed method, are important not only because of their applications in quantum theory of ferromagnetism~\cite{wmtb-ising-1976}, but also as a convenient benchmark for testing and comparing numerical algorithms and analytic approaches.
 Bailey, Borwein and Crandall~\cite{bailey-ising-2006} approached this problem from many different directions, and their results mark the state of the art of what can be achieved using the algorithms and methods of \nth{20} century.
 This is not an easy competition, and we are pleased that our algorithm stands up for it: we are able to reproduce the values calculated in~\cite{bailey-ising-2006} and also to improve the precision of physically relevant integrals from $5$--$6$ to $18$--$32$ decimal digits in dimensions $d\lesssim1000.$ 
 Using multiple precision library developed by David Bailey~\cite{bailey-mpfun}, we were able to reach precision of over $100$ decimal digits which revealed sub--exponential convergence of our algorithm $\eps\sim\exp(-\sqrt{\neval})$ for one of the considered integrals.
 The potential to converge sub-exponentially w.r.t. the number of function evaluations clearly distinguish the proposed method from MC/qMC algorithms, which usually demonstrate sublinear convergence $\eps\sim\neval^{-\gamma}$ with $0\leq\gamma\leq1.$
 
 The use of multiple precision calculations increases the challenge of high precision.
 Even though MPFUN2015~\cite{bailey-mpfun} and other arbitrary precision libraries~\cite{bailey-highprec-2015} are well optimised, the lack of optimisation at CPU level and vectorisation at the level of \texttt{BLAS} operations slows the calculations down, as well as requires extra steps when \texttt{BLAS} and \texttt{Lapack} functions need to be re-implemented in multiple precision.
 Although this problem is mitigated in more modern languages (such as Matlab, Python and Julia), they do not always provide enough control of parallelisation at both the distributed--memory (MPI) and shared--memory (OpenMP) levels.
 This is why for the development and demonstration stage we decided to implement the algorithm in Fortran, although it is clear that further work is required to simplify access to end users through interfaces to high--level languages mentioned above.

 The context of numerical integration is particularly convenient because the final answer is simply a number, allowing us to objectively evaluate and compare the quality of different algorithms for the given problem.
 It is good to see that for the examples we considered in this paper tensor cross interpolation is superior to MC and qMC algorithms.
 However it must be noted that the proposed method does not just compute the integral, but reconstructs the whole function in the high--dimensional tensor--product domain and represents it in TT form.
 When the compact representation of the function is available, it can be post--processed (e.g. interactively) to produce projections, nonlinear functionals (e.g. high--order moments), etc.
 This approach can be compared to calculation with functions using Chebyshev polynomials~\cite{trefethen-approx-2013}, and integrating Chebyshev interpolation together with the tensor cross interpolation seems to be a natural direction for further work, continuing the existing work in two and three--dimensions~\cite{trefethen-chebfun2-2013,trefethen-chebfun3-2017}.

 The most important direction of development of this work is without doubt the application of the proposed method to larger variety of applications.
 Many problems motivating precise high--dimensional integration are listed in~\cite{bailey-highprec-2015};
  we can extend this list by mentioning applications in 
    multivariate probability~\cite{dafs-tt-bayes-2018pre}, 
    stochastics~\cite{dklm-tt-pce-2015,ds-alscross-2019}, and
    optimal control~\cite{dp-chemotaxis-2018pre,qds-qcontrol-2019pre}.
 We are hopeful that the proposed tensor cross interpolation algorithm will demonstrate fast convergence in these applications and eventually becomes a method of choice for high--dimensional integration.


\section*{Software}
The Fortran implementation of Alg.~\ref{alg} is made by both authors and available at:
\begin{itemize}
 \item \href{http://github.com/savostyanov/ttcross}{\texttt{github.com/savostyanov/ttcross}}.
\end{itemize}


\begin{thebibliography}{10}

\bibitem{bailey-pslq-2000}
{\sc D.~H. Bailey}, \href {http://dx.doi.org/10.1109/5992.814653} {{\em Integer
  relation detection}}, Computing in Science \& Engineering, 2 (2000),
  pp.~24--28.

\bibitem{bailey-mpfun}
\leavevmode\vrule height 2pt depth -1.6pt width 23pt, \href
  {https://www.davidhbailey.com/dhbpapers/mpfun2015.pdf} {{\em {MPFUN2015:} {A}
  thread-safe arbitrary precision computation package}},  (2015).
\newblock https://www.davidhbailey.com/dhbpapers/mpfun2015.pdf.

\bibitem{bailey-highprec-2015}
{\sc D.~H. Bailey and J.~M. Borwein}, \href
  {http://dx.doi.org/10.3390/math3020337} {{\em High-precision arithmetic in
  mathematical physics}}, Mathematics, 3 (2015), pp.~337--367.

\bibitem{bailey-ising-2006}
{\sc D.~H. Bailey, J.~M. Borwein, and R.~E. Crandall}, \href
  {http://dx.doi.org/10.1088/0305-4470/39/40/001} {{\em Integrals of the
  {Ising} class}}, J Phys A: Math. Gen., 39 (2006), pp.~12271--12302.

\bibitem{bg-htcross-2015}
{\sc J.~Ballani and L.~Grasedyck}, {\em Hierarchical tensor approximation of
  output quantities of parameter-dependent {PDEs}}, {SIAM/ASA} Journal on
  Uncertainty Quantification, 3 (2015), pp.~852--872.

\bibitem{lars-htcross-2013}
{\sc J.~Ballani, L.~Grasedyck, and M.~Kluge}, \href
  {http://dx.doi.org/10.1016/j.laa.2011.08.010} {{\em Black box approximation
  of tensors in hierarchical {Tucker} format}}, Linear Algebra Appl., 428
  (2013), pp.~639--657.

\bibitem{lars-review-2014}
\leavevmode\vrule height 2pt depth -1.6pt width 23pt, \href
  {http://dx.doi.org/10.1007/978-3-319-08159-5_10} {{\em A review on adaptive
  low-rank approximation techniques in the hierarchical tensor format}}, in
  Extraction of Quantifiable Information from Complex Systems, vol.~102 of
  Lecture Notes in Computational Science and Engineering, Springer, 2014,
  pp.~195--210.

\bibitem{BarthSchwab-ml_mc_spde-2011}
{\sc A.~Barth, C.~Schwab, and N.~Zollinger}, \href
  {http://dx.doi.org/10.1007/s00211-011-0377-0} {{\em Multi-level {Monte Carlo}
  finite element method for elliptic {PDEs} with stochastic coefficients}},
  Numerische Mathematik, 119 (2011), pp.~123--161.

\bibitem{bartholdi-1982}
{\sc J.~J. Bartholdi~III}, {\em A good submatrix is hard to find}, Operations
  Research Lett., 1 (1982), pp.~190--193.

\bibitem{bebe-2000}
{\sc M.~Bebendorf}, \href {http://dx.doi.org/10.1007/pl00005410} {{\em
  Approximation of boundary element matrices}}, Numer. Math, 86 (2000),
  pp.~565--589.

\bibitem{bellman-dyn-program-1957}
{\sc R.~E. Bellman}, {\em Dynamic programming}, Princeton University Press,
  1957.

\bibitem{Bieri2009}
{\sc M.~Bieri and C.~Schwab}, \href
  {http://dx.doi.org/10.1016/j.cma.2008.08.019} {{\em {Sparse high order FEM
  for elliptic sPDEs}}}, Comp. Meth. Appl. Mech. Eng., 198 (2009),
  pp.~1149--1170.

\bibitem{mcmc-handbook-2011}
{\sc S.~Brooks, A.~Gelman, G.~Jones, and X.-L. Meng}, eds., {\em Handbook of
  {M}arkov chain {M}onte {C}arlo}, CRC Press, 2011.

\bibitem{griebel-sparsegrids-2004}
{\sc H.-J. Bungatrz and M.~Griebel}, \href
  {http://dx.doi.org/10.1017/S0962492904000182} {{\em Sparse grids}}, Acta
  Numerica, 13 (2004), pp.~147--269.

\bibitem{Schwab-HOQMC-2014}
{\sc J.~Dick, F.~Y. Kuo, Q.~T.~L. Gia, D.~Nuyens, and C.~Schwab}, \href
  {http://dx.doi.org/10.1137/130943984} {{\em Higher order {QMC}
  {P}etrov--{G}alerkin discretization for affine parametric operator equations
  with random field inputs}}, SIAM J. Num. An., 52 (2014), pp.~2676--2702.

\bibitem{Kuo-QMC-2013}
{\sc J.~Dick, F.~Y. Kuo, and I.~H. Sloan}, \href
  {http://dx.doi.org/10.1017/S0962492913000044} {{\em High-dimensional
  integration: The quasi-{M}onte {C}arlo way}}, Acta Numerica, 22 (2013),
  pp.~133--288.

\bibitem{dafs-tt-bayes-2018pre}
{\sc S.~Dolgov, K.~Anaya-Izquierdo, C.~Fox, and R.~Scheichl}, \href
  {http://arxiv.org/abs/1810.01212} {{\em Approximation and sampling of
  multivariate probability distributions in the tensor train decomposition}},
  {arXiv} preprint 1810.01212, 2018.

\bibitem{dkh-cme-2014}
{\sc S.~Dolgov and B.~Khoromskij}, \href {http://dx.doi.org/10.1002/nla.1942}
  {{\em Simultaneous state-time approximation of the chemical master equation
  using tensor product formats}}, Numer. Linear Algebra Appl., 22 (2015),
  pp.~197--219.

\bibitem{dklm-tt-pce-2015}
{\sc S.~Dolgov, B.~N. Khoromskij, A.~Litvinenko, and H.~G. Matthies}, \href
  {http://dx.doi.org/10.1137/140972536} {{\em {P}olynomial {C}haos {E}xpansion
  of random coefficients and the solution of stochastic partial differential
  equations in the {T}ensor {T}rain format}}, {SIAM} J. Uncertainty
  Quantification, 3 (2015), pp.~1109--1135.

\bibitem{dp-chemotaxis-2018pre}
{\sc S.~Dolgov and J.~W. Pearson}, \href {http://arxiv.org/abs/1806.08539}
  {{\em Preconditioners and tensor product solvers for optimal control problems
  from chemotaxis}}, {arXiv} preprint 1806.08539, 2018.

\bibitem{dpss-frac-2015}
{\sc S.~Dolgov, J.~W. Pearson, D.~V. Savostyanov, and M.~Stoll}, \href
  {http://dx.doi.org/10.1016/j.amc.2015.09.042} {{\em Fast tensor product
  solvers for optimization problems with fractional differential equations as
  constraints}}, Applied Mathematics and Computation, 273 (2016), pp.~604 --
  623.

\bibitem{ds-alscross-2019}
{\sc S.~Dolgov and R.~Scheichl}, \href {http://dx.doi.org/10.1137/17M1138881}
  {{\em A hybrid {A}lternating {L}east {S}quares--{TT-Cross} algorithm for
  parametric {PDEs}}}, SIAM/ASA Journal on Uncertainty Quantification, 7
  (2019), pp.~260--291.

\bibitem{DKhOs-parabolic1-2012}
{\sc S.~V. Dolgov, B.~N. Khoromskij, and I.~V. Oseledets}, \href
  {http://dx.doi.org/10.1137/120864210} {{\em Fast solution of
  multi-dimensional parabolic problems in the tensor train/quantized tensor
  train--format with initial application to the {Fokker}-{Planck} equation}},
  SIAM J. Sci. Comput., 34 (2012), pp.~A3016--A3038.

\bibitem{dks-ttfft-2012}
{\sc S.~V. Dolgov, B.~N. Khoromskij, and D.~V. Savostyanov}, \href
  {http://dx.doi.org/10.1007/s00041-012-9227-4} {{\em Superfast {Fourier}
  transform using {QTT} approximation}}, J. Fourier Anal. Appl., 18 (2012),
  pp.~915--953.

\bibitem{ds-dmrgamen-2015}
{\sc S.~V. Dolgov and D.~V. Savostyanov}, \href
  {http://dx.doi.org/10.1007/978-3-319-10705-9_33} {{\em Corrected one-site
  density matrix renormalization group and alternating minimal energy
  algorithm}}, in Numerical Mathematics and Advanced Applications --- {ENUMATH}
  2013, vol.~103, 2015, pp.~335--343.

\bibitem{dst-fb-2014}
{\sc S.~V. Dolgov, A.~P. Smirnov, and E.~E. Tyrtyshnikov}, \href
  {http://dx.doi.org/10.1016/j.jcp.2014.01.029} {{\em Low-rank approximation in
  the numerical modeling of the {F}arley-{B}uneman instability in ionospheric
  plasma}}, J. Comp. Phys., 263 (2014), pp.~268--282.

\bibitem{dongarra-top10}
{\sc J.~Dongarra and F.~Sullivan}, {\em Introduction to the top 10 algorithms},
  Computing in Science \& Engineering, 2 (2000), pp.~22--23.

\bibitem{mahoney-cur-2006}
{\sc P.~Drineas, R.~Kannan, and M.~W. Mahoney}, \href
  {http://dx.doi.org/10.1137/S0097539704442702} {{\em Fast {Monte Carlo}
  algorithms for matrices {III}: {Computing} a compressed approximate matrix
  decomposition}}, SIAM J Comput, 36 (2006), pp.~184--206.

\bibitem{Einstein-relativitat-1916}
{\sc A.~Einstein}, \href {http://dx.doi.org/10.1002/andp.19163540702} {{\em
  {D}ie {G}rundlage der allgemeinen {R}elativit\"{a}tstheorie}}, Annalen der
  Physik, 354 (1916), pp.~769--822.

\bibitem{etter-par-als-2016}
{\sc S.~Etter}, \href {http://dx.doi.org/10.1137/15M1038852} {{\em Parallel
  {ALS} algorithm for solving linear systems in the hierarchical tucker
  representation}}, SIAM J. Sci. Comput., 38 (2016), pp.~A2585--A2609.

\bibitem{fannes-mps-1992}
{\sc M.~Fannes, B.~Nachtergaele, and R.~Werner}, \href
  {http://dx.doi.org/10.1007/BF02099178} {{\em Finitely correlated states on
  quantum spin chains}}, Comm. Math. Phys., 144 (1992), pp.~443--490.

\bibitem{fkst-chem-2008}
{\sc H.-J. Flad, B.~N. Khoromskij, D.~V. Savostyanov, and E.~E. Tyrtyshnikov},
  \href {http://dx.doi.org/10.1515/RJNAMM.2008.020} {{\em Verification of the
  cross {3D} algorithm on quantum chemistry data}}, Rus. J. Numer. Anal. Math.
  Model., 23 (2008), pp.~329--344.

\bibitem{foster-1997}
{\sc L.~V. Foster}, \href {http://dx.doi.org/10.1016/S0377-0427(97)00154-4}
  {{\em The growth factor and efficiency of gaussian elimination with rook
  pivoting}}, J. Comput. Appl. Math., 86 (1997), pp.~177--194.

\bibitem{gantmacher-1959}
{\sc F.~R. Gantmacher}, {\em The theory of matrices}, Clelsea, NY, 1959.

\bibitem{golubkahan-1965}
{\sc G.~Golub and W.~Kahan}, {\em Calculating the singular values and
  pseudo-inverse of a matrix}, SIAM J. Numer. Anal., 2 (1965), pp.~205--224.

\bibitem{golub-2013}
{\sc G.~Golub and C.~Van~Loan}, {\em Matrix computations}, Johns Hopkins
  University Press, Baltimore, MD, 2013.

\bibitem{gostz-maxvol-2010}
{\sc S.~A. Goreinov, I.~V. Oseledets, D.~V. Savostyanov, E.~E. Tyrtyshnikov,
  and N.~L. Zamarashkin}, {\em How to find a good submatrix}, in Matrix
  Methods: Theory, Algorithms, Applications, V.~Olshevsky and E.~Tyrtyshnikov,
  eds., World Scientific, Hackensack, NY, 2010, pp.~247--256.

\bibitem{gt-skel-2011}
{\sc S.~A. Goreinov and E.~E. Tyrtyshnikov}, \href
  {http://dx.doi.org/10.1134/S1064562411030355} {{\em Quasioptimality of
  skeleton approximation of a matrix in the {Chebyshev} norm}}, Doklady Math.,
  83 (2011), pp.~374--375.

\bibitem{gt-psa-1995}
{\sc S.~A. Goreinov, E.~E. Tyrtyshnikov, and N.~L. Zamarashkin}, {\em
  Pseudo--skeleton approximations of matrices}, Reports of Russian Academy of
  Sciences, 342 (1995), pp.~151--152.

\bibitem{gtz-psa-1997}
\leavevmode\vrule height 2pt depth -1.6pt width 23pt, \href
  {http://dx.doi.org/10.1016/S0024-3795(96)00301-1} {{\em A theory of
  pseudo--skeleton approximations}}, Linear Algebra Appl., 261 (1997),
  pp.~1--21.

\bibitem{gtz-maxvol-1997}
{\sc S.~A. Goreinov, N.~L. Zamarashkin, and E.~E. Tyrtyshnikov}, \href
  {http://dx.doi.org/10.1007/BF02358985} {{\em Pseudo--skeleton approximations
  by matrices of maximum volume}}, Mathematical Notes, 62 (1997), pp.~515--519.

\bibitem{graham-QMC-2011}
{\sc I.~Graham, F.~Kuo, D.~Nuyens, R.~Scheichl, and I.~Sloan}, \href
  {http://dx.doi.org/10.1016/j.jcp.2011.01.023} {{\em {Q}uasi-{M}onte {C}arlo
  methods for elliptic {PDEs} with random coefficients and applications}}, J.
  Comput. Phys., 230 (2011), pp.~3668--3694.

\bibitem{gras-hsvd-2010}
{\sc L.~Grasedyck}, \href {http://dx.doi.org/10.1137/090764189} {{\em
  Hierarchical singular value decomposition of tensors}}, SIAM J. Matrix Anal.
  Appl., 31 (2010), pp.~2029--2054.

\bibitem{grasedyck-par-cross-2015}
{\sc L.~Grasedyck, R.~Kriemann, C.~L{\"o}bbert, A.~N{\"a}gel, G.~Wittum, and
  K.~Xylouris}, \href {http://dx.doi.org/10.1007/s00791-015-0247-x} {{\em
  Parallel tensor sampling in the hierarchical {Tucker} format}}, Computing and
  Visualization in Science, 17 (2015), pp.~67--78.

\bibitem{gu-rrqr-1995}
{\sc M.~Gu and C.~Eisenstat}, \href {http://dx.doi.org/10.1137/0917055} {{\em
  Efficient algorithms for computing a strong rank--revealing {QR}
  factorization}}, SIAM J. Sci. Comput., 17 (1996), pp.~848--869.

\bibitem{hackbusch-2012}
{\sc W.~Hackbusch}, {\em Tensor Spaces And Numerical Tensor Calculus},
  Springer--Verlag, Berlin, 2012.

\bibitem{hk-ht-2009}
{\sc W.~Hackbusch and S.~{K\"uhn}}, \href
  {http://dx.doi.org/10.1007/s00041-009-9094-9} {{\em A new scheme for the
  tensor representation}}, J. Fourier Anal. Appl., 15 (2009), pp.~706--722.

\bibitem{trefethen-chebfun3-2017}
{\sc B.~Hashemi and L.~N. Trefethen}, \href
  {http://dx.doi.org/10.1137/16M1083803} {{\em Chebfun in three dimensions}},
  SIAM J. Sci. Comput., 39, p.~C341–C363.

\bibitem{bokh-surv-2015}
{\sc B.~N. Khoromskij}, \href {http://dx.doi.org/10.1051/proc/201448001} {{\em
  Tensor numerical methods for multidimensional {PDE}s: theoretical analysis
  and initial applications}}, ESAIM: Proc., 48 (2015), pp.~1--28.

\bibitem{klumper-mps-1993}
{\sc A.~Kl\"umper, A.~Schadschneider, and J.~Zittartz}, \href
  {http://dx.doi.org/10.1209/0295-5075/24/4/010} {{\em Matrix product ground
  states for one-dimensional spin-1 quantum antiferromagnets}}, Europhys.
  Lett., 24 (1993), pp.~293--297.

\bibitem{knuth-1985}
{\sc D.~E. Knuth}, \href {http://dx.doi.org/10.1080/03081088508817636} {{\em
  Semi--optimal bases for linear dependencies}}, Linear and Multilinear
  Algebra, 17 (1985), pp.~1--4.

\bibitem{kolda-review-2009}
{\sc T.~G. Kolda and B.~W. Bader}, \href {http://dx.doi.org/10.1137/07070111X}
  {{\em Tensor decompositions and applications}}, SIAM Rev., 51 (2009),
  pp.~455--500.

\bibitem{kumar-review-2017}
{\sc K.~N. Kumar and J.~Schneider}, \href
  {http://dx.doi.org/10.1080/03081087.2016.1267104} {{\em Literature survey on
  low rank approximation of matrices}}, Linear and Multilinear Algebra, 65
  (2017), pp.~2212--2244.

\bibitem{Scheichl-mlqmc-lognorm-2017}
{\sc F.~Kuo, R.~Scheichl, C.~Schwab, I.~Sloan, and E.~Ullmann}, \href
  {http://dx.doi.org/10.1090/mcom/3207} {{\em Multilevel quasi-{M}onte {C}arlo
  methods for lognormal diffusion problems}}, Math. Comp.,  (2017),
  pp.~2827--2860.

\bibitem{Schwab-MLQMC-2015}
{\sc F.~Y. Kuo, C.~Schwab, and I.~H. Sloan}, \href
  {http://dx.doi.org/10.1007/s10208-014-9237-5} {{\em Multi-level quasi-monte
  carlo finite element methods for a class of elliptic pdes with random
  coefficients}}, Found. Comp. Math.,  (2015), pp.~1--39.

\bibitem{metropolis-mc-1949}
{\sc N.~Metropolis and S.~Ulam}, \href
  {http://dx.doi.org/10.1080/01621459.1949.10483310} {{\em The {Monte} {Carlo}
  method}}, Journal of the American statistical association, 44 (1949),
  pp.~335--341.

\bibitem{meyer-mctdh-1990}
{\sc H.-D. Meyer, U.~Manthe, and L.~S. Cederbaum}, {\em The
  multi-configurational time-dependent {H}artree approach}, Chem. Phys. Lett.,
  165 (1990), pp.~73--78.

\bibitem{mo-rectmaxvol-2018}
{\sc A.~Y. Mikhalev and I.~V. Oseledets}, \href
  {http://dx.doi.org/10.1016/j.laa.2017.10.014} {{\em Rectangular
  maximum-volume submatrices and their applications}}, Linear Algebra Appl.,
  (2018), pp.~187--211.

\bibitem{morokoff-qmc-1995}
{\sc W.~J. Morokoff and R.~E. Caflisch}, \href
  {http://dx.doi.org/10.1006/jcph.1995.1209} {{\em Quasi-{Monte} {Carlo}
  integration}}, J Comp. Phys., 122 (1995), pp.~218--230.

\bibitem{poole-1992}
{\sc L.~Neal and G.~Poole}, \href
  {http://dx.doi.org/10.1016/0024-3795(92)90432-A} {{\em A geometric analysis
  of {Gaussian} elimination. {II}}}, Linear Alg. Appl., 173 (1992),
  pp.~239--264.

\bibitem{nieder-qmc-1978}
{\sc H.~Niederreiter}, {\em Quasi--{Monte Carlo} methods and pseudo--random
  numbers}, Bull. AMS, 84 (1978), pp.~957--1041.

\bibitem{nobile-sg-mc-2016}
{\sc F.~Nobile, L.~Tamellini, F.~Tesei, and R.~Tempone}, {\em An adaptive
  sparse grid algorithm for elliptic {PDE}s with lognormal diffusion
  coefficient}, in Sparse Grids and Applications - Stuttgart 2014, Springer
  International Publishing, 2016, pp.~191--220.

\bibitem{osel-tt-2011}
{\sc I.~V. Oseledets}, \href {http://dx.doi.org/10.1137/090752286} {{\em
  Tensor-train decomposition}}, SIAM J. Sci. Comput., 33 (2011),
  pp.~2295--2317.

\bibitem{DoOs-dmrg-solve-2011}
{\sc I.~V. Oseledets and S.~V. Dolgov}, \href
  {http://dx.doi.org/10.1137/110833142} {{\em Solution of linear systems and
  matrix inversion in the {TT}-format}}, SIAM J. Sci. Comput., 34 (2012),
  pp.~A2718--A2739.

\bibitem{ost-tucker-2008}
{\sc I.~V. Oseledets, D.~V. Savostianov, and E.~E. Tyrtyshnikov}, \href
  {http://dx.doi.org/10.1137/060655894} {{\em Tucker dimensionality reduction
  of three-dimensional arrays in linear time}}, SIAM J. Matrix Anal. Appl., 30
  (2008), pp.~939--956.

\bibitem{ost-chem-2010}
{\sc I.~V. Oseledets, D.~V. Savostyanov, and E.~E. Tyrtyshnikov}, \href
  {http://dx.doi.org/10.1002/nla.682} {{\em Cross approximation in tensor
  electron density computations}}, Numer. Linear Algebra Appl., 17 (2010),
  pp.~935--952.

\bibitem{ot-ttcross-2010}
{\sc I.~V. Oseledets and E.~E. Tyrtyshnikov}, \href
  {http://dx.doi.org/10.1016/j.laa.2009.07.024} {{\em {TT-cross} approximation
  for multidimensional arrays}}, Linear Algebra Appl., 432 (2010), pp.~70--88.

\bibitem{pan-rrlu-2000}
{\sc C.-T. Pan}, \href {http://dx.doi.org/10.1016/S0024-3795(00)00120-8} {{\em
  On the existence and computation of rank--revealing {LU} factorizations}},
  Linear Algebra Appl., 316 (2000), pp.~199--222.

\bibitem{poole-2000}
{\sc G.~Poole and L.~Neal}, \href
  {http://dx.doi.org/10.1016/S0377-0427(00)00406-4} {{\em The rook's pivoting
  strategy}}, J. Comput. Appl. Math., 123 (2000), pp.~353--369.

\bibitem{qds-qcontrol-2019pre}
{\sc D.~Qui\~nones Valles, S.~Dolgov, and D.~Savostyanov}, \href
  {http://arxiv.org/abs/1903.00064} {{\em Tensor product approach to quantum
  control}}, {arXiv} preprint 1903.00064, 2019.

\bibitem{rst-volterra-2014}
{\sc J.~A. Roberts, D.~V. Savostyanov, and E.~E. Tyrtyshnikov}, \href
  {http://dx.doi.org/10.1016/j.cam.2013.10.025} {{\em Superfast solution of
  linear convolutional {Volterra} equations using {QTT} approximation}}, J.
  Comput. Appl. Math., 260 (2014), pp.~434--448.

\bibitem{sav-rr-2009}
{\sc D.~V. Savostyanov}, \href {http://dx.doi.org/10.4208/nmtma.2009.m9006s}
  {{\em Fast revealing of mode ranks of tensor in canonical form}}, Numer.
  Math. Theor. Meth. Appl., 2 (2009), pp.~439--444.

\bibitem{sav-rank1-2012}
\leavevmode\vrule height 2pt depth -1.6pt width 23pt, \href
  {http://dx.doi.org/10.1016/j.laa.2011.11.008} {{\em {QTT}-rank-one vectors
  with {QTT}-rank-one and full-rank {Fourier} images}}, Linear Algebra Appl.,
  436 (2012), pp.~3215--3224.

\bibitem{sav-qott-2014}
\leavevmode\vrule height 2pt depth -1.6pt width 23pt, \href
  {http://dx.doi.org/10.1016/j.laa.2014.06.006} {{\em Quasioptimality of
  maximum--volume cross interpolation of tensors}}, Linear Algebra Appl., 458
  (2014), pp.~217--244.

\bibitem{sdwk-nmr-2014}
{\sc D.~V. Savostyanov, S.~V. Dolgov, J.~M. Werner, and I.~Kuprov}, \href
  {http://dx.doi.org/10.1103/PhysRevB.90.085139} {{\em Exact {NMR} simulation
  of protein-size spin systems using tensor train formalism}}, Phys. Rev. B, 90
  (2014), p.~085139.

\bibitem{so-dmrgi-2011proc}
{\sc D.~V. Savostyanov and I.~V. Oseledets}, \href
  {http://dx.doi.org/10.1109/nDS.2011.6076873} {{\em Fast adaptive
  interpolation of multi-dimensional arrays in tensor train format}}, in
  Proceedings of 7th International Workshop on Multidimensional Systems (nDS),
  IEEE, 2011.

\bibitem{schneider-cross2d-2010}
{\sc J.~Schneider}, \href {http://dx.doi.org/10.1016/j.jat.2010.04.012} {{\em
  Error estimates for two--dimensional cross approximation}}, J. Approx.
  Theory, 162 (2010), pp.~1685--1700.

\bibitem{smolyak-1963}
{\sc S.~A. Smolyak}, {\em Quadrature and interpolation formulas for tensor
  products of certain class of functions}, Dokl. Akad. Nauk SSSR, 148 (1963),
  pp.~1042--1053.
\newblock Transl.: Soviet Math. Dokl. 4:240-243, 1963.

\bibitem{solomonik-cyclops-2014}
{\sc E.~Solomonik, D.~Matthews, J.~R. Hammond, J.~F. Stanton, and J.~Demmel},
  \href {http://dx.doi.org/10.1016/j.jpdc.2014.06.002} {{\em A massively
  parallel tensor contraction framework for coupled-cluster computations}},
  Journal of Parallel and Distributed Computing, 74 (2014), pp.~3176 -- 3190.

\bibitem{stewart-lr-2000}
{\sc G.~W. Stewart}, \href {http://dx.doi.org/10.1109/5992.814658} {{\em The
  decompositional approach to matrix computation}}, Computing in Science \&
  Engineering, 2 (2000), pp.~50--59.

\bibitem{white-parallel-dmrg-2013}
{\sc E.~Stoudenmire and S.~White}, \href
  {http://dx.doi.org/10.1103/PhysRevB.87.155137} {{\em Real-space parallel
  density matrix renormalization group}}, Phys. Rev. B, 87 (2013), p.~155137.

\bibitem{stuart-bayes-2010}
{\sc A.~M. Stuart}, \href {http://dx.doi.org/10.1017/S0962492910000061} {{\em
  Inverse problems: A {B}ayesian perspective}}, Acta Numerica, 19 (2010),
  pp.~451--559.

\bibitem{Tadmor-exp_acc_diff-1986}
{\sc E.~Tadmor}, {\em The exponential accuracy of {Fourier} and {Chebychev}
  differencing methods}, SIAM J. Numer. Anal., 23 (1986), pp.~1--23.

\bibitem{trefethen-chebfun2-2013}
{\sc A.~Townsend and L.~N. Trefethen}, \href
  {http://dx.doi.org/10.1137/130908002} {{\em An extension of {Chebfun} to two
  dimensions}}, SIAM J. Sci. Comput., 35 (2013), p.~C495–C518.

\bibitem{trefethen-approx-2013}
{\sc L.~N. Trefethen}, {\em Approximation Theory and Approximation Practice},
  SIAM, 2013.

\bibitem{tee-mosaic-1996}
{\sc E.~E. Tyrtyshnikov}, \href {http://dx.doi.org/10.1007/BF02575706} {{\em
  Mosaic-skeleton approximations}}, Calcolo, 33 (1996), pp.~47--57.

\bibitem{tee-cross-2000}
\leavevmode\vrule height 2pt depth -1.6pt width 23pt, \href
  {http://dx.doi.org/10.1007/s006070070031} {{\em Incomplete cross
  approximation in the mosaic--skeleton method}}, Computing, 64 (2000),
  pp.~367--380.

\bibitem{white-dmrg-1992}
{\sc S.~R. White}, \href {http://dx.doi.org/10.1103/PhysRevLett.69.2863} {{\em
  Density matrix formulation for quantum renormalization groups}}, Phys. Rev.
  Lett., 69 (1992), pp.~2863--2866.

\bibitem{wilkinson-1961}
{\sc J.~H. Wilkinson}, \href {http://dx.doi.org/10.1145/321075.321076} {{\em
  Error analysis of direct method of matrix inversion}}, J Assoc. Comp.
  Machinery, 8 (1961), pp.~281--330.

\bibitem{wmtb-ising-1976}
{\sc T.~T. Wu, B.~M. McCoy, C.~A. Tracy, and E.~Barouch}, \href
  {http://dx.doi.org/10.1103/PhysRevB.13.316} {{\em Spin-spin correlation
  functions for the two-dimensional {Ising} model: {Exact} theory in the
  scaling region}}, Phys. Rev. B, 13 (1976), p.~316.

\bibitem{yang-ising-1952}
{\sc C.~N. Yang}, \href {http://dx.doi.org/10.1103/PhysRev.85.808} {{\em The
  spontaneous magnetization of a two-dimensional {Ising} model}}, Phys. Rev.,
  85 (1952), p.~808.

\bibitem{zo-maxvol-2017}
{\sc N.~L. Zamarashkin and A.~I. Osinsky}, \href
  {http://dx.doi.org/10.1134/S1064562416060156} {{\em New accuracy estimates
  for pseudoskeleton approximations of matrices}}, Doklady Mathematics, 94
  (2016), pp.~643--645.

\bibitem{zyokd-ttcircuit-2015}
{\sc Z.~Zheng, X.~Yang, I.~V. Oseledets, G.~E. Karniadakis, and L.~Daniel},
  \href {http://dx.doi.org/10.1109/TCAD.2014.2369505} {{\em Enabling
  high-dimensional hierarchical uncertainty quantification by {ANOVA} and
  {Tensor-Train} decomposition}}, IEEE Trans. Comput-aided Des. Integr.
  Circuits Syst., 34 (2015), pp.~63--76.

\end{thebibliography}

\end{document}